\documentclass[11pt]{article}

\usepackage{amssymb}
\usepackage{latexsym}
\usepackage{xypic}
\usepackage[mathscr]{euscript}
\usepackage{graphicx}
\usepackage{tikz-cd}
\usepackage[T1]{fontenc}
\usepackage[utf8]{inputenc}
\usepackage{suetterl}
\usepackage{amsthm}
\usepackage{bm}

\usepackage{amscd}
\usepackage{amsmath}
\usepackage{amssymb}
\usepackage{fullpage}

\usepackage{amssymb}
\usepackage{latexsym}
\usepackage[mathscr]{euscript}
\usepackage{latexsym}
\usepackage{amsfonts}
\usepackage{epsfig}
\usepackage{enumitem}

\newcommand\mb{\mathbb} 
\newcommand\mc{\mathcal}

\def \bG{\mathbf G}
\def \rmG{\rm G}

\def \bK{\mathbf K}
\def \O{\mathcal O}

\def \rmX{\rm X}

\newcommand\commentout[1]{\marginpar{\tiny $\backslash$commentout}}

\newcommand\mfk{\mathbf G}

\setlist[enumerate]{label=(\arabic*)}

\newtheorem{Lemma}{Lemma}[section]
\newtheorem{Theorem}{Theorem}[section]
\newtheorem{Proposition}{Proposition}[section]
\newtheorem{Definition}{Definition}[section]
\newtheorem{Corollary}{Corollary}[section]

\newtheorem{Example}{Example}[section]
\newtheorem{Remark}{Remark}[section]

\newenvironment{Proof}{\par\noindent\textbf{Proof:}}
{\qed}

\usepackage{graphics}

\title{Representation Theory and the Algebraic $K$-Theory of Geometric Fields}
\author{Gunnar Carlsson \footnote{Research supported in part by NSF DMS-0406992}  \\ Department of Mathematics, Stanford University  \\Stanford, California 94305 \and 
Roy Joshua \footnote{Research supported in part by Simons Foundation MP-TSM-00002471}\\ Department of Mathematics, The Ohio State University \\ Columbus, Ohio 43210}

\begin{document}
\maketitle
\tableofcontents

\section{Introduction}
One of the central themes in number theory and field theory is the relationship between arithmetic properties of a field and properties of Galois groups of extensions of the field.  Obvious examples include Kummer theory and class field theory, but of course the proven  Bloch-Kato and Milnor conjectures  and the entire Langlands program fall under this umbrella.        

In this paper we will  study another relationship of this type, which will relate the representation theory of the absolute Galois group $G_F$ of a geometric field $F$ (i.e. $F$ containing an algebraically closed subfield) to its higher algebraic $K$-theory $K_*(F)$.  To describe this relationship, we will need to introduce a derived version of the completion of the representation ring $R[G_F]$ at the kernel of its mod $l$ augmentation. 

Let $f:A \rightarrow B$ be a homomorphism of commutative rings, and let $M$ be an $A$-module. 
$M^{\wedge}_B$.  The output of this construction is a simplicial $A$-module, and one may study its homotopy groups.  For the case of Noetherian rings $A$, surjective ring homomorphisms $f$ and finitely generated $A$-modules $M$, this construction yields the usual notion of  completion in the sense that $\pi _0M_B$ is isomorphic to the $I_f$-adic completion of $M$, where $I_f$ denotes the kernel of $f$, and where all the higher homotopy groups vanish. 

 When $A$ is non-Noetherian, however, $\pi _* M^{\wedge}_B$ is often  non-zero for $*>0$. In  particular, for profinite groups $G$, we define the representation ring of $R[G]$ to be the colimit of the representation rings of the finite quotients of $G$ defining the profinite structure, and find that for non-finite profinite groups, the representations rings are usually non-Noetherian. This fact allows for derived completions of these rings to have non-trivial higher homotopy groups.  What is very surprising is that these  homotopy groups have very interesting geometric interpretations when the group $G$ is the pro-$l$ completion of a discrete group $\Gamma$, where $l$ denotes a prime. 
 
 In \cite{lawson}, it is shown that for a nilpotent group $\Gamma$,  the homology  groups of  a spectrum constructed from a stabilized version of the representation varieties of $\Gamma$ called {\em deformation $K$-theory} are  identical to the homotopy groups of the completion $R[\hat{\Gamma}_{\ell}]_{\mathbb{F}_{\ell}}^{\wedge}$, where $\Gamma _{\ell}$ denotes the pro-$l$ completion of $\Gamma$.  For example for the group $\mathbb{Z}^n$, the homotopy groups of the derived completion of the representation ring of $\mathbb{Z}_{\ell}^n$ are isomorphic to an exterior algebra over $\mathbb{Z}_{\ell}$ on $n$ one-dimensional generators.  Similarly, for the free pro-$l$ group on $n$ generators, we obtain that $\pi _0 \cong \mathbb{Z}_{\ell}$, $\pi _1 \cong \mathbb{Z}_{\ell}^n$, and $\pi _i \cong 0$ for $i \geq 2$. The idea that  purely discrete algebraic constructions on representation rings of finite groups allows one to reconstruct continuous representation spaces of infinite discrete groups  is by itself  extremely  surprising. The further finding  that  these constructions applied to absolute or separable Galois groups also yield arithmetic information is even more striking.  

 Recall the definition of {\em Milnor $K$-theory}  $K_*^M(F) $ of a field $F$ from \cite{milnor}. Note that it depends only on the arithmetic in the field, with no reference to Galois actions on ambient fields.  
  With these definitions, we are able to state the following. 
\begin{Theorem}\label{bigtheorem} 
    Let $F$ denote a  field  of characteristic zero.  Suppose further that $F$ contains an algebraically closed subfield $k$. Let $G_F$ denote the absolute   Galois group of $F$.   Let $R[G_F]$ denote the ring of isomorphism classes of finite dimensional $k$-linear representations of $G_F$.   Then we have  isomorphisms 
    $$ \pi _i ({R}[G_F]^{\wedge}_{\mb{F}_{\ell}}) \cong K^{M}_i(F) ^{\wedge}_{\ell}
    $$ where $K^M_*$ denotes Milnor $K$-theory,  $l$ is a prime, and where ${R}[G_F]^{\wedge}_{\mb{F}_{\ell}}$ denotes the derived completion of $R[G]$ at the mod l augmentation homomorphism $R[G] \rightarrow \mb{F}_{\ell}$.   The isomorphism is functorial with respect to inclusions of subfields.  
\end{Theorem}

For $i =1$, this statement is equivalent to Kummer theory.  The arithmetic interpretations of the higher homotopy groups  of $R[G_F]^{\wedge}_{\mb{F}_{\ell}}$ should be thought of as a higher dimensional and non-commutative extension of Kummer theory.

\begin{Remark}
    This result extends in a straightforward way to fields of characteristic $p$ which admit no $p$-power extensions, and to separable Galois groups, when $\ell \neq p$.  For clarity of exposition, we only treat the characteristic zero case.  
\end{Remark}

This theorem about homotopy groups follows from the following spectrum level theorem.  
In \cite{completion}, it is shown  that the derived completion construction extends to $E_{\infty}$ ring spectra and  modules over them,  as well as  homomorphisms of them.  

We will adopt the following terminology for the remainder of the paper. $k$ will always denote an algebraically closed field of characteristic zero.  For $G$ a finite group and $X$ a scheme of finite type over $k$, we will denote by $\bK(X,G) $ and $\mfk(X,G)$ the equivariant  $K$- and $G$-theory spectra of $X$ as defined in \cite[ section 2]{royg2}. See Section \ref{equivariantk} for a discussion of properties we will need for these constructions. An extension of these constructions to profinite group actions and pro-$G$-schemes is carried out in Section \ref{extensions}.  $\bK(X,G) $ is an $E_{\infty}$ ring spectrum, and $\mfk (X,G) $ is a module spectrum over it.

\begin{Theorem}\label{nextbigtheorem} 
Let the notation and hypotheses of Theorem \ref{bigtheorem} hold. Then we have the following.  
 \begin{enumerate}[label=(\arabic*), font=\upshape] \item{The map of spectra 
 $$\mfk(Spec(k), G_F )_{\mb{F}_{\ell} }^{\wedge} \rightarrow \mfk(Spec(E), G_F )^{\wedge} _{\mb{F}_{\ell}} 
 $$
 is an equivalence, where $E$ denotes the algebraic  closure of $F$,   and where $\mb{F}_{\ell}$ denotes the Eilenberg-MacLane spectrum for the finite field $\mb{F}_{\ell}$. 
 Moreover,
 $$\mfk(Spec(k), G_F) \simeq \bK(Spec(k), G_F) \mbox{ and } \mfk(Spec(E), G_F) \simeq \bK(Spec(E), G_F).$$}
 \item{There is a natural  equivalence  $\bK(Spec(F))^{\wedge}_{l} \rightarrow \bK(Spec(E), G_F )^{\wedge} _{\mb{F}_{\ell}}   $, where $\bK(Spec(F)) ^{\wedge}_{\ell}$ denotes the Bousfield-Kan $l$-adic completion (see \cite{bousfield}) of the $K$-theory spectrum of $F$.  It follows that the $l$-completed  $K$-groups of $F$ are computable as the homotopy groups of the derived completion of the $\bK$-theory of  the spectrum of the category of finite dimensional  $k$-linear representations of $G_F$, at its augmentation to $\mb{F}_{\ell}$.
 }
 \end{enumerate}
 
\end{Theorem}

Theorem \ref{nextbigtheorem} is the key technical result which allows us to prove the algebraic result in Theorem \ref{bigtheorem}.

{Algebraic topology solves geometric problems  using algebraic invariants, and spectral sequences are an important tool for computing those algebraic invariants.  In actually solving geometric problems, computing algebraic invariants is only part of the story.  Typically, it is also important to have further information about how the algebraic invariants related to geometric properties of the  objects at hand, and they are often not captured entirely by the algebraic structure of $E_2$-terms of  spectral sequences.  For instance, the behavior of symmetries of a geometric object are often an important ingredient.  }

{The revolutionary work of Voevodsky on the conjectures of Milnor and Bloch-Kato (see \cite{voevodsky} and \cite{weibel}) identifies the $E_2$-term of  the Bloch-Lichtenbaum spectral sequence converging to the algebraic $K$-theory of a field.  However,   it does not provide a geometric model that explicitly incorporates the behavior of the absolute Galois group on algebraic closure.    Our construction is built directly from only the separable  Galois group $G$ and the $K$-theory of a separably closed field with $G$-action, i.e the spectrum $ku^{\wedge}_{\ell}$. We believe strongly in the potential for algebraic $K$-theory to solve arithmetic and geometric problems, and that the existence of a model for it constructed directly from $G$ will be critical  in this effort.    }

The authors wish to thank A. Blumberg,  B. Conrad, P. Diaconis, D. Dugger, L. Hesselholt, T. Lawson, G. Lyo,  I. Madsen, D. Ramras, R. Vakil, and K. Wickelgren for numerous helpful conversations on this work.  We especially thank Ben Mann for his support and extensive input on this project.

\section{Motivation/Discussion} Since the results of \cite{completion} have not been broadly used to study the $K$-theory of fields, we think it will be useful to clarify the relationship between our results and the standard descent approaches based on homotopy fixed point spectra.  
\subsection{Generalities on homotopy fixed point sets}
We elaborate on the  second point above.
To understand the relationship of our result with the Quillen-Lichtenbaum conjecture, we'll describe some computations of homotopy fixed point sets in two other situations.  

We recall that given a group $G$ and an action of $G$ on a space $X$, we define the homotopy fixed point set  (or space or spectrum)  of the action is the function space of equivariant maps $F^G(EG, X)$ from $EG$, a contractible space on which $G$ acts freely, to $X$, which we denote by $X^{hG}$.  The evident projection $EG \rightarrow *$ induces a map $ \epsilon: X^G \rightarrow X^{hG}$, since  equivariant maps of $*$ to $X$ are simply fixed points of the $G$-action. 

An extremely useful feature of homotopy fixed point sets is that in the case of an action  based at a point $x_0 \in X$,  their homotopy groups can be computed via a spectral sequence with $E_2$ term $H^{-p}(G,\pi _q (X,x_0))$, where $\pi _q(X,x_0)$ is a $G$-module via the induced action of $G$.  Note that the $E_2$-term does not require any knowledge of the homotopy type of the fixed point set $X^G$, but only of the group action on the homotopy groups of the ambient space $X$.  The map $\epsilon$ is rarely an equivalence, but it has been studied in a number of situations.  The first instance occurred in the work of Atiyah and Segal \cite{atiyahsegal}, where the homotopy fixed point spectrum of the equivariant periodic $K$-theory  spectrum $K_G(*)$ for a finite group $G$ was completely analyzed. The equivariant $K$-theory spectrum has homotopy groups $\pi _*KU(*) \otimes R[G]$, where $R[G]$ denotes the complex representation ring of $G$, and it was proved in \cite{atiyahsegal} that the homotopy groups of the homotopy fixed point spectrum are given by $$\pi _* KU(*) \otimes R[G]^{\wedge}_{\rho}$$
where $R[G]^{\wedge}_{\rho} $ denotes the completion of $R[G]$ at its augmentation ideal.  
  The conclusion is that although the map $\epsilon$ is not an equivalence, it is in a sense close to one, and it differs from one in a  predictable way.  

A corresponding result for stable homotopy theory  based on a conjecture of G.B. Segal is found in \cite{segalconjecture}, where the finding is again that the homotopy groups of the homotopy fixed point set on the action on equivariant stable cohomotopy theory $\pi _G^*(*)$ for a finite group are obtained as completions at the augmentation ideals of the {\em Burnside ring}, over which the groups $\pi_G^*(*) $ form a graded algebra. 

There is also an unstable result for finite $G$-complexes ($G$ again finite)  called the {\em Sullivan conjecture}, which can be found in  \cite{sullivanconj}, \cite{sullivanconj2}. The moral from all these results is that for each case, the fixed point set can be recovered with small and predictable inaccuracies from the homotopy fixed point set, and that those inaccuracies are  ``erased" by a well understood completion process.
\subsection{The special case of algebraic $K$-theory}

The idea behind the Quillen-Lichtenbaum conjecture was that this should also be true when $X$ is the algebraic $K$-theory spectrum of the algebraic closure of a field $F$, equipped with the action of $G_F$, the absolute Galois group of $F$.  However, it was almost immediately observed that this naive idea failed as stated, and the conjecture was formulated as a statement which held in sufficiently high dimensions. In order to address the failure of the Quillen-Lichtenbaum conjecture to provide a precise  computation of the $K$-theory of fields, we have to come to grips with it  and   modify  the conjecture.

 The modification  comes through the invention of {\em motivic homotopy theory} by F. Morel and V. Voevodsky \cite{motivichomotopy}.  In addition to developing analogues of the ideas of standard homotopy theory which instead applies to categories of schemes, the authors of \cite{motivichomotopy} construct two contractible objects on which an algebraic group acts.  One of these is the straightforward analogue of the usual bar construction  $EG$ in the non-motivic situation, but the other, denoted $E_{gm}G$, the ``geometric universal space for $G$",  has no counterpart in the standard setting. One can then define the ``geometric homotopy fixed  point set" of the algebraic $K$-theory $KX$ of a variety or scheme $X$ with action by an algebraic group  to be the space of equivariant maps from $E_{gm}G$ to $KX$, and ask whether the natural map from $(KX)^G$ to this homotopy fixed point set is better behaved. 
 
 It turns out that it is, as is shown for algebraic group actions in \cite{royg2}, where analogues of the results of Atiyah and Segal \cite{atiyahsegal} are proved for the algebraic $K$ theory of schemes with action by an algebraic group $G$. What was shown  in \cite{royg4} is that these results also extend to case of profinite groups and the action of the profinite absolute Galois group $G_F$ of a geometric field.  In short, that the problem with the original statement of the Quillen-Lichtenbaum conjecture was that the homotopy fixed points  were computed in the {\em wrong category}, i.e. in the standard category of spaces rather than in the motivic category.

\subsection{Completion of spectra versus completion of homotopy groups}

 We assume we are given a field $F$ containing an algebraically closed subfield $k \subseteq F$, let $\overline{F} $ denote its algebraic closure, and let $G_F$ denote its absolute Galois group $Gal(\overline{F} /F)$, in general a profinite group.  In Section \ref{extensions} we construct a version of the equivariant $K$-theory spectrum,  $\bK(-,G_F)$,  which applies to profinite groups and certain actions on pro-schemes, and observe that $\bK(Spec (\overline{F}), G_F) \cong \bK(Spec (F))$.  One may also consider the natural map of ring spectra $\alpha_F : \bK(Spec(k), G_F) \rightarrow 
 \bK(Spec (\overline{F}), G_F)$ induced by the inclusion $k \hookrightarrow \overline{F} $. 
 
 This map is far from being an equivalence since, for example, $\pi _0 \bK(Spec(k), G_F) $ is the free abelian group on the continuous $k$-linear representations of $G_F$, where $k$ is given the discrete topology, and on the other hand 
$$ \pi _0 \bK(Spec(\overline{F}), G_F) \cong \pi _0 \bK(Spec(F)) \cong {\mathbb Z}
$$
However, the spectra $\bK(Spec(k), G_F)$ and $\bK(Spec (\overline{F}), G_F)$ both carry multiplicative structures, which makes $\bG(Spec (\overline{F}), G_F)$ into a module over the commutative ring spectrum $\bK(Spec(k), G_F)$.  See \cite{ekmm} and \cite{symmetric} for a discussion of these multiplicative theories. 

We now apply the derived completion completion construction for commutative ring spectra introduced in \cite{completion} to the following situation.  We will denote by $\rho $ the homomorphism $\bK(Spec(k),G_F) \rightarrow \mb{F}_{\ell} $, where $\mb{F}_{\ell}$ denotes the mod $\ell$ Eilenberg-MacLane spectrum.  $\rho$ is induced by the functor which assigns to every representation of $G_F$ its dimension mod $l$, where $\mb{F}_{\ell} $ is regarded as a symmetric monoidal category with object set equal to $\mb{F}_{\ell} $, and with only identity morphisms.

We  construct derived completions along the homomorphism $\rho $, to obtain a map 
$$\alpha _F^{\wedge} : \bK(Spec(k), G_F)^{\wedge}_{\rho} \rightarrow \bK(Spec(\overline{F}), G_F)^{\wedge}_{\rho }$$
called the {\em representational assembly}, and the conjecture is that this map is an equivalence of spectra. In this paper, we prove this conjecture, and the resulting Theorem \ref{nextbigtheorem} is what allows us to prove Theorem \ref{bigtheorem}.

\section{Equivariant Algebraic  K-theory } \label{equivariantk}
Throughout the paper,  $k$ will denote an algebraically closed ground field.  We will consider the Waldhausen category model  for the algebraic $K$-theory of the categories of complexes of vector bundles  on schemes (see \cite{waldhausen} and \cite{TT}), which is in a natural way a symmetric spectrum (see \cite{symmetric}, \cite{geisser}, and \cite{blumberg}, Section 2). Specifically, we deal with the category of pseudo-coherent complexes (perfect complexes) of $\O$-modules provided with an action by the given
group $\rmG$ so that the differentials of the complex are equivariant and with bounded coherent cohomology
sheaves to define the $G$-equivariant $G$-theory spectrum $\bG(\quad, G)$ (the $\rmG$-equivariant $K$-theory spectrum $\bK(\quad, G)$.    In \cite{blumberg}, Corollary 2.7, conditions are identified that guarantee that algebraic $K$-theory spectra are $E_{\infty}$-ring spectra and modules over such, and in \cite{zakharevich}, it is shown that these conditions are satisfied.    For finite groups $G$  we  will be considering $G$-equivariant versions $\mfk(X, G) $, where the category associated to an affine  scheme $X = Spec(A_X)$  with $G$-action will consist of finitely generated $A_X$-modules $M$ with a $G$-action satisfying $ (\alpha \cdot m)^g = \alpha^g \cdot m^g$. For a closed $G$-subscheme $Y \subseteq X$, we will write
$\mfk_Y(X, G)$ for the $G$-equivariant $\mfk$-theory of the category of modules supported on $Y$.

\begin{Proposition}\label{kproperties} We enumerate the properties of $K$-theory spectra that we will need.
\begin{enumerate}[label=(\arabic*), font=\upshape]
    \item{{\bf Multiplicative structure:} \label{relequiv}$\bK(X, G)$ is an $E_{\infty}$ ring spectrum, with the multiplication constructed from tensor products of modules over $X$.  Moreover, $\bK_Y(X, G) $ and $\bG_{Y}(\rmX, G)$  are  $\bK(X, G)$-module spectra.   $\pi _0 \bK(k, G) \cong R[G]$, where $R[G]$ denotes the $k$-linear representation ring of $G$. There is a natural equivalence  of $\bK(X, G)$-modules $i_Y^X: \mfk(Y, G) \rightarrow \mfk_Y(X, G)$. If we have $Y \subseteq X^{\prime} \subseteq X$, then the natural map $\mfk_Y(X, G) \rightarrow \mfk_Y(X^{\prime}, G)$ is an equivalence. }
    \item{{\bf Poincar\'e-Duality for regular Noetherian schemes}: The obvious forgetful functor sending a perfect complex
    of equivariant $\O_X$-modules to the underlying pseudo-coherent complex induces a map of spectra:
    $\bK(X, G) \rightarrow \mfk(X, G)$. This map is a weak-equivalence when $X$ is a regular Noetherian scheme:
    in fact, the same proof in \cite[3.21 Theorem]{TT} for the non-equivariant case works here as well.}
    \item{{\bf Functoriality in $X$:} The $\mfk(-, G)$-construction is contravariantly functorial for equivariant flat maps of schemes, inducing homomorphisms of $E_{\infty}$-module spectrum maps.   }
    \item{{\bf Functoriality in $Y$:} $\mfk_Y(X, G)$ is covariantly functorial in $Y$, so $Y \subseteq Y^{\prime} \subseteq X$ induces a  map of $\bK(X, G)$-module-spectra  $\mfk_Y(X, G) \rightarrow  \mfk_{Y^{\prime}}(X, G)$.  }
    \item{{\bf Localization sequence:} \label{excisionintro} For a pair of $k$-schemes with $G$ action $(X,Y)$, where $Y$ is closed in $\rmX$, there is a fibration sequence of $\bK(X, G)$-module spectra 
    $$ \mfk_Y(X, G) \rightarrow \mfk(X , G) \rightarrow \mfk(X - Y, G) 
    $$ More generally, if we have $Z \subseteq Y \subseteq X$, with $Y$ and $Z$ closed in $X$, there is a fibration sequence 
    $$ \bG_Z(X, G) \rightarrow \bG_Y (X, G) \rightarrow \bG_{Y-Z} (X-Z, G).
    $$}
    \item{{\bf Homotopy property:} For a linear action of $G$ on $\mb{A}^n$, the natural map $\mfk(X, G) \rightarrow \mfk(X \times \mb{A}^n, G)$ induced by the projection $X \times \mb{A}^n\rightarrow X  $ is an equivalence of spectra.  }
    
    \item{{\bf Functoriality in $G$:} For a homomorphism of groups $G \rightarrow G^{\prime}$ and a $G^{\prime}$ action on a scheme $X$, there is an induced map of $\bK(Spec(k),G^{\prime})$-module spectra $\mfk(X, G^{\prime}) \rightarrow \mfk(X, G)$, where in the right hand term the action of $G$ on $X$ is obtained by pulling back the $G^{\prime}$-action.  } 
    \item{{\bf Eilenberg-MacLane spectra:} We let $\mb{F}_{\ell}$ denote the model for the mod-$\ell$ Eilenberg-MacLane spectrum as a symmetric spectrum.  It is an $E_{\infty}$-ring spectrum.   $\bK (Spec(k), \{e\} ) $ is equipped with a unique morphism $\varepsilon _{\ell}:\bK(Spec(k), \{*\}) \rightarrow \mb{F}_{\ell}$ of $E_{\infty}$-ring spectra. For any group $G$, we define the augmentation $\mfk(Spec(k), G) \rightarrow \mb{F}_{\ell}$ to be the composite $\bK(Spec(k), G) \rightarrow \bK(Spec(k), \{ e \}) \xrightarrow{\varepsilon _{\ell}} \mb{F}_{\ell}$. It induces the mod $l$ augmentation $R[G] \rightarrow \mb{F}_{\ell}$  on $\pi _0$. (We will allow ourselves to use the notation $\mb{F}_{\ell}$ both for the Eilenberg-MacLane spectrum for $\mb{F}_{\ell}$ and  for the ordinary ring $\mb{F}_{\ell}$. Which one  is meant will always be clear from the context.) }
    \item{{\bf Derived constructions:}  Because the category of symmetric spectra  admits a strictly commutative smash  product, it is possible to construct derived smash products of $\bK(Spec(k), G)$ modules $M$ and $N$. This is done by restricting the ``positive  model structure" on symmetric spectra to the category of  $\bK(Spec(k), G)$-module-spectra, and performing cofibrant replacement both of $\bK(Spec(k))$-algebras and modules over them.  We will treat every construction as taking place on cofibrant objects without further mention, since no confusion will result.    There results  a spectral sequence with $E_2$-term 
    $$Tor_*^{\bK_*(Spec(k), G)}(\pi _*M,\pi _*N) 
    $$
    converging to $\pi _* (M  \wedge _{\bK(Spec(k), G)} N )$. }
\end{enumerate}
    
\end{Proposition}
\section{Extensions to Profinite Groups and Pro-Schemes}\label{extensions}

We will need to extend the definition of $\mfk(X, G)$ to  profinite groups and pro-schemes.  In the case of profinite groups acting on schemes through a finite quotient, there is a very simple extension of equivariant $G$-theory. For $G$ a profinite group $G$, we say an action of $G$ on a scheme $X$ is {\em uniformly continuous} or {\em $u$-continuous} if the $G$-action factors through a finite quotient of $G$.  Letting $G$ be a profinite group, with u-continuous  action of $G$ on a scheme $X$, we simply define $\bK(X, G)$ ($\mfk(X, G)$) to be the colimit of the spectra $\bK(X, G/H)$ ($\mfk(X, G/H)$, respectively), as $H$ varies over all closed normal subgroups of finite index in $K$.  $\bK(Spec(k), G)$ becomes a commutative ring spectrum, and $\mfk(X, G)$ becomes a $\bK(Spec(k), G)$-module spectrum.  We refer to this as the {\em simple extension} to profinite groups. 
\begin{Remark}
    Strictly speaking,  one must replace the category of schemes $X$ with a $G$-action which factors through $G/K$, where $K$ is a closed normal subgroup of finite index with the following ``blown up'' category $\mc{S}^G$.  The objects of $\mc{S}^G$ are pairs  $(K,X)$, where $K$ is a closed subgroup of $G$ of finite index.  A   morphism from $(K,X) $ to $(K^{\prime},X^{\prime}) $ consists of a pair $(i,f)$, where $i : K \hookrightarrow K^{\prime}$ is an inclusion of subgroups of $G$ and where $f: X \rightarrow X^{\prime}$ is a $G$-equivariant map of schemes when $G$ acts on $X$ and $X^{\prime}$ via the projections $G \rightarrow G/K$ and $G \rightarrow G/K^{\prime}$. For simplicity of exposition, we will ignore this distinction, since its inclusion does not affect any of the arguments.    
\end{Remark}
We will need something more general, where we incorporate actions on pro-schemes which do not necessarily stabilize at any particular finite index subgroup.   We will  not require the full generality of all possible actions in the pro-category, but work instead in the following situation.  
\begin{Definition} \label{procontinuous}
    Let $G$ be a profinite group.  By a {\em pro-continuous $G$-scheme}, we mean a pro-object in the category of $u$-continuous $G$-schemes.   For a pro-continuous $G$-scheme $\tilde{X}$,  we define $\bK(\tilde{X},G)$ and $\mfk(\tilde{X}, G)$ to be the colimits $\underset{i \in I}{colim} \hspace{.05cm} \bK(\tilde{X}_i,G) $ and $\underset{i \in I}{colim}  \hspace{.05cm} \mfk(\tilde{X}_i,G) $, where $I$ is the cofiltered category parametrizing $\tilde{X}$.  $\bK(\tilde{X},G)$ is a commutative ring spectrum, and $\mfk(\tilde{X},G)$ is a $\bK(\tilde{X},G)$-module spectrum.  $Spec (k)$ denotes the constant $G$-scheme with value $Spec(k)$ and equipped with trivial $G$-action.  
\end{Definition}

Let  $k$ denote a field.    We make the following definitions in order to have a mechanism to construct the pro-continuous $G$-spectra we need.  
\begin{Definition}
   Let $A$ be a Noetherian $k$-algebra, and let $B$ be a $A$-algebra.  We say an $A$-subalgebra of $B^{\prime} \subseteq B $ is {\em strongly flat under $B$} if {\em(a)} it is faithfully flat and finite over $A$ and {\em (b)} $B$ is a faithfully flat $B^{\prime}$-algebra.  We say $B$ is a {\em reachable} $A$-algebra if it is the colimit of its strongly flat $A$-subalgebras. We denote by $SF(B)$ the category of all strongly flat $A$-subalgebras of $B$.  $SF(B)$ is a filtering category.  
\end{Definition}

With this definition it is clear that given a reachable $A$-algebra $B$, we obtain a pro-object in the category $Sch_k$ of schemes over $Spec(k)$, parametrized by $SF(B)^{op}$.    In order to apply $\bK$- and $\mfk$-theory, though, we will need to show that this pro-object actually has image in the category $Sch_k^{flat}$, i.e. the category of all schemes over $Spec(k)$ and flat morphisms.  The following proposition verifies  this. 

\begin{Proposition}\label{Bfunctoriality} Suppose that we have a reachable $A$-algebra $B$, and that $B_1 \subseteq B_2 $ is an inclusion of strongly flat  $A$-subalgebras under  $B$.  Then $B_2$ is a flat $B_1$-algebra.  It follows that  for any morphism $f:B^{\prime} \rightarrow B^{\prime \prime}$ of $SF(B)$, the corresponding morphism $Spec(f): Spec(B^{\prime \prime}) \rightarrow Spec(B^{\prime})$ is faithfully flat.  
    \end{Proposition}
    \begin{Proof} We need to prove that if $f: M \rightarrow N$ is an injective map of $B_1$-modules, then $id_{B_2} \otimes _{B_1} f$ is an injective map  of $B_2$-modules. But because $B$ is a flat $B_1$-module, we know that $id_{B} \otimes_{B_1}  f$ is injective.  But we have 
    $$ id_{B} \otimes_{B_1}  f  = id_B \otimes _{B_2}(id_{B_2} \otimes _{B_1} f)
    $$
    so because of the faithful flatness of $B$ as a $B_2$-module, $id_{B_2} \otimes _{B_1} f$ is injective. Faithful flatness is immediate because of the surjectivity of $Spec(B_2) \rightarrow Spec(B_1)$, which follows from the surjectivity of $Spec(B) \rightarrow Spec(B_1)$
    \end{Proof}

\begin{Definition}
    Let $A$ be a Noetherian $k$-algebra over a field $k$, and $B$  a flat   $A$-algebra, equipped with an action by a profinite group $G$  under which $A$ is invariant and the action of $G$ on  $Spec(A)$ is u-continuous.        
    \begin{enumerate}[label=(\arabic*), font=\upshape]
     \item{We say the $G$-action on $B$ is {\em  continuous } if for every element $b \in B$, the stabilizer of $b$ is a closed subgroup of finite index.}
    \item{We say $B$ is {\em $G$-reachable} if it is the colimit of its  strongly flat $G$-invariant $A$-subalgebras.   }
    \end{enumerate}
   In this situation we define $SF(B,G)$ to be the  category of strongly flat $G$-invariant $A$-subalgebras of $B$.  It is a filtering category. 
\end{Definition}
\begin{Proposition} \label{actiontrivial}
    Given $A$ and $B$ as in the definition above, and let $B^{\prime}$ be a strongly flat $G$-invariant $A$-subalgebra  of $B$.  Then $Spec(B^{\prime})$ is u-continuous.    
\end{Proposition}
\begin{Proof}
    $B^{\prime}$ is a finitely generated algebra over $A$, and we have assumed that the $G$-action on $Spec(A)$ is u-continuous. Let $\{ b_1, \ldots , b_n \}$ be a generating set for $B^{\prime}$ as an $A$-algebra.    The continuity of the $G$-action on $B$ then shows that the subgroup of $G$ which stabilizes $B^{\prime}$ is the intersection of the stabilizer of  $A$ with the  intersection of all the subgroups $Stab(b_i )$, each of which is closed and of finite index.  The result follows.  
\end{Proof}

We will also have occasion to use functoriality in the algebra $A$.
\begin{Definition} \label{specialty}
       Let $B$ be an $A$-algebra.  Then we say that $(A,B)$ is a {\em $G$-reachable pair}  if $B$ is a $G$-reachable $A$-algebra.  By a morphism of $G$-reachable pairs from $(A,B) $ to $(A^{\prime}, B^{\prime})$, we mean a pair of  faithfully flat morphisms $f_A: A \rightarrow A^{\prime}$ and $f_B: B \rightarrow B^{\prime} $ so that $i^{\prime} \circ  f_A = f_B \circ  i$, where $i:A \rightarrow B$ and $i^{\prime}: A ^{\prime} \rightarrow B^{\prime}$ are the inclusions. Associated to  a morphism $(f_A, f_B): (A,B) \rightarrow (A^{\prime}, B^{\prime}) $ is an inclusion $i_{(f_A, f_B)}: A^{\prime}   \otimes _{A} B\rightarrow B^{\prime}$ of $A^{\prime}$-algebras. 
       
        \end{Definition}

        \begin{Proposition}\label{funccriterion} 
Let $f: (A,B)\rightarrow (A^{\prime },B^{\prime}) $ be a map of $G$-reachable pairs.   Then for any strongly flat $A$-subalgebra $C$  under $B$, the $A^{\prime}$-subalgebra $f(C) $ is strongly flat under $B^{\prime} $
  
\end{Proposition}
\begin{Proof}
   This is an elementary exercise in the behavior of faithful flatness under pushouts and compositions. \end{Proof}

   \begin{Example}
    Let $F$ be a geometric field containing an algebraically closed subfield $k$ , and let $\overline{F}$ denote its algebraic closure.  Let $G_F$ denote the absolute Galois group of $F$.  Then $(F, \overline{F})$ is a $G_F$-reachable pair.  
\end{Example}
\begin{Example}
    Suppose $A$ is a Noetherian $k$-algebra equipped with a u-continuous action by a profinite group $G$.  Then $(A,A) $ is a $G$-reachable pair. 
\end{Example}
\begin{Example} The Puiseux algebras defined and used in  Sections  \ref{Calculation} and \ref{Bextensions} are $G$-reachable algebras, when equipped with a Puiseux action as defined in Section \ref{Calculation}.  
\end{Example} 

Our goal is to attach  a pro-continuous $G$-scheme to a $G$-reachable $A$-algebra $B$, to which we can then apply $\mfk$-theory as in Definition \ref{procontinuous}. 

\begin{Definition}\label{bigdefinition} Let $(A,B)$ be a $G$-reachable pair, where $G$ is a profinite group.  
The category $SF(B,G)$ is filtering, and it follows that the contravariant functor $B^{\prime} \rightarrow Spec(B^{\prime})$ is a pro-object in the category of schemes with $G$-action. It follows from Proposition \ref{Bfunctoriality}  that it has image in the subcategory of flat maps of $G$-schemes, and from Proposition \ref{actiontrivial} that it is pro-continuous.  We will denote this pro-$G$-scheme by $Spec(A,B) $, and by $\mfk(Spec(A,B ),G)$ the corresponding $\bK(Spec(k),G)$-module spectrum $\mfk (Spec(A,B), G)$, as in Definition \ref{procontinuous}.  $Spec(-)$ is a contravariant functor from the category of $G$-reachable pairs to the category of pro-continuous $G$-schemes.  
\end{Definition}
\begin{Remark}
    The notation $Spec$ may appear a bit overloaded, but it can be viewed as an extension of $Spec(-)$ from the category of u-continuous Noetherian affine  $G$-schemes $X$ embedded as the $G$-reachable pairs of the form $(\mc{A}(X), \mc{A}(X))$, where $\mc{A}$ denotes the global sections functor on $G$-schemes, to the category of all $G$-reachable pairs. We will sometimes write $Spec(B)$ for $Spec(A,B)$, when $A$  is understood.  In particular, when $\overline{F}$ is the algebraic closure of $F$, with $G_F$ the absolute Galois group, $Spec(\overline{F})$ will denote the pro-$G_F$ scheme of all the finite Galois extensions of $F$.  
\end{Remark}

\begin{Proposition} \label{simpleproperties} $G$-reachability enjoys the following two properties. 
   \begin{enumerate}
       
\item{Let $B_1$ and $B_2$ be algebras over $A_1$ and $A_2$, and suppose that both $B_1$ and $B_2$ are equipped with continuous $G$-actions, restricting to the identity on $A_1$ and $A_2$.   Suppose that both $B_1$ and $B_2$ are $G$-reachable.  Then $B_1 \otimes _k B_2$ is a $G$-reachable $A_1 \otimes _k A_2$-algebra. On pro-$G$-schemes, we will write $Spec(B_1 \otimes _k B_2) \cong Spec(B_1) \times _k Spec(B_2) $.  }
\item{Let $B$ be an $A$-algebra, and $G$ a profinite group.  Suppose $B$ is equipped with a continuous $G$ action, and that $A$ is invariant under the action and is fixed elementwise by a closed subgroup $H \subseteq G$  of finite index.   Suppose also that $B$ is $G$-reachable.   Let $G^{\prime}\rightarrow   G$ be a continuous homomorphism of profinite groups.  Then $B$ is $G^{\prime}$-reachable.  }
   \end{enumerate}
\end{Proposition}
\begin{Proof}
    Elementary. 
\end{Proof}

Given a $G$-reachable $A$-algebra $B$, and a closed subscheme $Y \subseteq Spec(A)$, we will now construct two additional module spectra over $\bK(k,G)$.   

\begin{Definition} \label{pullbackschemes}
    Let $A$ be a Noetherian $k$-algebra, and let $B$ be an $A$-algebra with continuous $G$-action leaving $A$ invariant, and so that there is a finite index closed subgroup $H \subseteq G $ fixing $A$ elementwise.  Suppose further that $B$ is  $G$-reachable. Let $\pi :Spec(A,B) \rightarrow Spec(A) $ denote the projection of pro-continuous $G$-schemes, where $Spec(A)$ is regarded as a constant pro-scheme.   For any $G$-invariant subscheme $Y \subseteq Spec(A)$, we let $Spec(A,B; Y)$ denote the pro-continuous $G$-scheme given by the functor on $SF(B,G) $ defined  by $B^{\prime} \rightarrow \pi _{B^{\prime}}^{-1}(Y)$, where $\pi _{B^{\prime}}: Spec(B^{\prime}) \rightarrow Spec(A)$ is the projection. If the subscheme $Y$ is affine with affine coordinate ring $A_Y$,  then $B \otimes _{A} A_Y$ is itself a $G$-reachable  $A_Y$-algebra, and $Spec^{}(A,B; Y) \cong Spec(A_Y, B \otimes  A_Y)$.  When $Y$ is closed, we may define  the $\bK(Spec(k),G)$-module spectrum $\mfk _{Y} (Spec(A, B),G)$ to be the colimit of the  Ind-spectrum $B^{\prime} \rightarrow \mfk_{\pi ^{-1}_{B^{\prime}}Y}(Spec(B^{\prime}), G)$. 
    \end{Definition}

\begin{Proposition} \label{pullbacklocalization}
When $Y \subseteq Spec(A)$ is a $G$-invariant closed subscheme,  we obtain a fibration  sequence of $\bK (Spec(k),G)$-module spectra 
$$ \mfk_{Y}(Spec(A, B),G) \rightarrow \mfk (Spec(A, B),G) \rightarrow \mfk(Spec(A,B; Spec(A) - Y),G) 
$$
and more generally, when $Z \subseteq Y$ is an inclusion of $G$-invariant closed subschemes of $Spec(A)$, we obtain the localization sequence 
$$ \mfk_{Z}(Spec(A, B),G) \rightarrow \mfk_Y(Spec(A, B),G) \rightarrow \mfk_{Y-Z}(Spec (A, B ; Spec(A)-Z),G) 
$$
\end{Proposition}
\begin{Proof}
    This is a straightforward verificaton using item \ref{excisionintro} of Proposition \ref{kproperties}. 
\end{Proof}

   \begin{Corollary}\label{pairfunctoriality}
       For a morphism $(A,B) \rightarrow (A^{\prime},B^{\prime})$ of $G$-reachable pairs, there is a natural map of $G$-pro-schemes $Spec(A^{\prime}, B^{\prime}) \rightarrow Spec(A, B)$ and consequently of the corresponding $G$-theory spectra, turning the $Spec$-construction into a contravariant functor on the category of $G$-reachable pairs.  It follows that the corresponding $G$-theory constructions are given a covariant $\bK(Spec(k),G) $-module spectrum valued functor on the category of $G$-reachable pairs.  
   \end{Corollary}

We conclude this section with a technical result. $\bK(Spec(k), G)$ is a commutative ring spectrum, and we consider the  commutative $\bK(Spec(k), G)$-algebras $C$, in particular $C = \mb{F}_{\ell}$ and $C = \bK(Spec(k), \{ e \} )$.  
We suppose that $\tilde{X}$ is a pro-continuous $G$-scheme, and can therefore be described as a $G$-scheme valued functor $I \rightarrow Sch_k$ for a cofiltering category.  We apply $\mfk(-;G)$, and observe that $\mfk(\tilde{X},G)$ is a colimit over spectra $\mfk(X(i), G(i))$ where $i$ is an object of $I$, and $G(i)$ is a finite quotient of $G$.  
\begin{Lemma}\label{tensorcolim}
    The spectrum $C \wedge_{\mfk(Spec(k),G)}\mfk(\tilde{X},G)$ is equivalent to the colimit $$\underset{i \in I}{colim}\hspace{.05cm}  \left[ C \wedge _{\mfk(Spec(k),G(i))} \mfk(X(i), G(i)) \right]$$ 
\end{Lemma}
    \begin{Proof}
        Immediate from standard colimit properties of derived tensor products.  
    \end{Proof}

\section{G-theory of Fields}

We will also need to understand in more detail the Galois-equivariant $\mfk$-theory of fields.  Let $E/F$ be a Galois extension of fields (possibly infinite), with Galois group $G$ (in general profinite).  As observed in Section \ref{extensions}, $E$ is a $G$-reachable $F$-algebra with continuous $G$-action.  $\bK(E, G) $ is constructed from the Waldhausen category  of ``descent data", i.e. $E$-vector spaces $V$ equipped with an  $F$-linear $G$-action $\rho$  for which $\rho(g) (e \cdot v)  = e^g \cdot \rho(g)(v)$. If $F$ contains  $k$ as a subfield,  there is a functor  $\iota $ from the category of finite dimensional $k$-linear representations $Rep_k [G]$ to $\mathcal{D}(E,G)$ given by $\iota(V,\rho) =  E \otimes _kV$, equipped with the $G$-action $\rho _E(g) (e \otimes v ) = e^g \otimes \rho (g)v$. It is immediate that the $\bK(\quad, G)$-theory spectrum of $k $ is the $\bK$-theory of the category of finite dimensional $k$-linear $G$-representations, that tensor product of representations,  makes $\bK(Spec(k), G) $ into a commutative ring spectrum.  The   functor $\iota$ makes  $\bK(Spec(E), G)$ into a module spectrum over $\bK(Spec(k), G)$. That module spectrum structure is induced by the pairing which assigns to a pair $(\sigma , (V, \rho ))$ the descent datum $(\sigma \otimes _k V, \sigma \otimes \rho)$.  

\begin{Proposition}\label{descentproperties}  Let $k$ be algebraically closed.  Then these constructions have the following properties. 
\begin{enumerate}[label=(\arabic*), font=\upshape]
    \item{ $\pi _0 \bK(Spec(k), G)   \cong R[G]$, where $R[G]$ denotes the  representation ring of $G$ on $k$-vector spaces. If $k$ is of characteristic zero, then it is also isomorphic to the usual complex representation ring.  }  
    \item{\label{fixedfield} $\bK(Spec(E), G) \cong KF$, so $K_*F$ becomes an $R[G]$-module via the   map $\iota$.  }
    \item{\label{zeroness} Let $I \subseteq R[G]$ denote the augmentation ideal.  Then $I \cdot K_*F = \{ 0 \} $}
\end{enumerate}
    
\end{Proposition}
\begin{Proof}
    Items (1) and (2)  are standard.  For item (3), any finite dimensional $k$-representation $(W,\sigma)  $  of $G$, as an element of $R[G] \cong \pi _0\bK(Spec(k), G)$, produces the  functor $T_{\sigma }: \mathcal{D}(E,G) \rightarrow \mathcal{D}(E,G) $ which is defined by 
    $$ T_{\sigma }(V,\rho) = (W \otimes _k V, \sigma \otimes \rho) 
    $$
 It is standard that the $\mfk$-theory of the category $\mathcal{D}(E,G)$ is naturally equivalent to the $\mfk$-theory of the field $F$, and that this equivalence is realized by a the functor $j$ which assigns to an $F$-vector space $U$ the descent data $(E \otimes _F U)$, equipped with the action $g \cdot (e \otimes u) = (g\cdot e) \otimes u$.  Since $j$ induces an equivalence of spectra, if we can show that $T_{\sigma} \circ j = T_{\sigma ^0} \circ j$, where $\sigma ^0$ denotes the underlying space of $\sigma$ equipped with the trivial $G$-action, we will have our result.  But $(T_{\sigma} \circ j)(U)  = (\sigma \otimes E) \otimes _F U$.  The result now follows from the standard fact that $\sigma \otimes _k E \cong \sigma ^0 \otimes E$. 
\end{Proof}
\section{Derived Completions} \label{derivedsection} 
  We recall the results of \cite{completion} concerning the definition and properties of the derived completion construction.  Given a homomorphism $ A \rightarrow B$ of commutative ring spectra and a $A$-module $M$,  a derived completion construction $M^{\wedge}_B$, which is itself a module spectrum over $A$.  If $M$ is actually an $A$-module, then so is $M^{\wedge} _B$. 

 \begin{Proposition} \label{properties} The derived completion construction enjoys the following properties.  
     
\begin{enumerate}[label=(\arabic*), font=\upshape] 
\item{\label{tripledefinition} $M^{\wedge} _B$ is defined as the total spectrum of the cosimplicial spectrum of the monad $M \rightarrow T_A(M;B) =   B \wedge_A M$ on the category of $A$-modules.  We write $\mc{T}_A^. (-;B)$ for the associated cosimplicial spectrum, and $\tilde{\mc{T}}_A( -;B)$ for the pro-spectrum $\{ Tot^i(\mc{T}_A^. (-;B))\}_i $. We extend $\tilde{\mc{T}}_A$ to a monad on the category of pro-module spectra in the evident way.  }
     \item{ \label{keyprop} For a homomorphism of $A$-modules $f:M \rightarrow N$, if the induced homomorphism of $A$-modules $B \wedge _A M \overset{id_B \wedge f}{\longrightarrow} B \wedge_A N$ is an equivalence then 
     $$M^{\wedge} _B \overset{\hat{f}_B}\longrightarrow N^{\wedge}_B
     $$ is an equivalence.  }
     \item{For actual Noetherian commutative rings $A$ and $B$ (discrete commutative ring spectra) and finitely generated  $A$-modules $M$, and with  $f$  surjective with kernel $I$, then $M^{\wedge} _B$ is also discrete, and is the usual completion $M^{\wedge} _I$ at the ideal $I$.   }
     \item{\label{derss}For a homomorphism $f:A \rightarrow B$ of commutative ring spectra and an $A$-module spectrum $M$, there is a second quadrant spectral sequence with $E_2$-term 
     $$ E_1^{pq} = \pi  _{q + 2p}(( \pi _{ -p}M)^{\wedge}_{\pi _0 B})
$$ converging to $\pi _{p+q}(M^{\wedge}_B)$}
\item{\label{trivialequiv}  Let $f: A \rightarrow B$ be a surjective map of commutative ring spectra. Let  $M$ be  a $B$-module.   Then the natural homomorphism  $M \rightarrow M^{\wedge}_B$ is an equivalence.  }
     \end{enumerate}
     \end{Proposition}

We will need to understand the functoriality properties of the derived completion construction in the map $ A \rightarrow B$.  In order to do this, we will introduce the notation $M^{\wedge}_{A \rightarrow B}$ for the completion construction associated to the morphism - it is identical to the construction $M^{\wedge}_B$ above when the ring $A$ is understood.   

\begin{Proposition} \label{moreproperties} The construction $M \rightarrow M^{\wedge}_{A \rightarrow B}$ satisfies the following functoriality properties.  
\begin{enumerate} [label=(\arabic*), font=\upshape]
\item{\label{primaryexact} The construction $M \rightarrow M^{\wedge}_{A \rightarrow B} $ is functorial for maps of $A$-modules. It preserves cofibration sequences of $A$-modules, as well as finite colimits.   }
    \item{\label{functormodule} Given maps of  of commutative ring spectra $A \rightarrow B \rightarrow C $,  there is an induced  morphism of $A$-modules $M^{\wedge}_{A \rightarrow C} \rightarrow M^{\wedge}_{B\rightarrow C} $.
    }
    \item{The completion $A^{\wedge}_{A \rightarrow C}$ is a commutative ring spectrum, and the composite map $A^{\wedge}_{A \rightarrow C} \rightarrow B^{\wedge}_{A \rightarrow C} \rightarrow B^{\wedge}_{B \rightarrow C} $ is a map of commutative ring spectra.  }
    \item{For a commutative ring spectrum $A$ and an $A$-module $M$, $M^{\wedge} _{id: A \xrightarrow{}A}\cong M $. }

\end{enumerate}
    
\end{Proposition}
\begin{Proof}
    These are all contained in Proposition 3.2 of \cite{completion}. 
\end{Proof}

 \begin{Lemma}\label{sphereequiv} Let $f: A \rightarrow \mb{F}_{\ell} $ be map of commutative ring spectra, and let $I$ denote the kernel of $\pi _0(f)$.  Let $M$ be an $A$-module spectrum so that $I$ acts trivially on $\pi _*(\mb{F}_{\ell} \wedge _{S^0} M)$.  Then the natural map $M^{\wedge}_{S^0 \rightarrow \mb{F}_{\ell}}  \rightarrow M^{\wedge}_{A \rightarrow \mb{F}_{\ell}} $ is an equivalence of spectra.  
 \end{Lemma} 
 \begin{Proof}
     From the spectral sequence of item \ref{derss} of Proposition \ref{properties}, it will suffice to prove that the map 
 $$ \pi_i (M)^{\wedge}_{\pi _0(S^0) \rightarrow \mb{F}_{\ell} } \rightarrow \pi _i(M)^{\wedge} _{\pi _0(A) \rightarrow \mb{F}_{\ell} }$$ is an equivalence for all $i$. We construct the monads $T_1= \mc{T}_{\pi _0(A)}(-; \mb{F}_{\ell}) $ and $T_2= \mc{T}_{\pi _0(S^0) }(-; \mb{F}_{\ell}) = \mc{T}_{\mb{Z}}(-;\mb{F}_{\ell} )$ on the category of $\pi _0(A)$-modules. Of course $\pi _0(A)$, $\pi _0(S^0)$, and $\pi _i(M)$ are all actual rings and modules, but as such we regard them as  ring spectra and module spectra in order to  use them as input to derived completion. There is a natural transformation $T_2 \rightarrow T_1$.   We  construct the bicosimplicial spectrum $X^{\cdot \cdot}$ defined by $X^{p,q}  \rightarrow T_1^{p+1}T_2 ^{q+1}(M) $.  The total spectrum  of $X^{..}$ is naturally equivalent to the total spectrum of $T^{\cdot}_1(M)$ by Lemma 2.14 of \cite{completion}, since there is a natural transformation of monads $T_2 \rightarrow T_1.  $ On the other hand, fixing $q$, we obtain we obtain the cosimplicial spectrum $T_1^{.} T_2^{q+1}(M) $, whose total spectrum is the completion of $T_2^{q+1}(M)$.  By the hypothesis on $M$, though, it follows that the map 
 $$ T_2^{q+1}(M) \rightarrow \underset{p}{Tot}\hspace{.05cm}T_1 ^{p+1}T_2 ^{q+1}(M) 
 $$
 is an equivalence using item \ref{trivialequiv} of Proposition \ref{properties}, since the $A$-module structure on $T^{q+1}(M)$ factors through the ring homomorphism $A \rightarrow \mb{F}_{\ell}$.  The result now follows.  
 \end{Proof}

 \begin{Corollary}\label{finalcriterion} Let $F$ be a geometric field, $\overline{F}$ its algebraic closure, and $G$ its absolute Galois group.  Then the natural map of $\bK(Spec(k))^{\wedge} _{\mb{F}_{\ell}}$-modules  
      $$\mfk(Spec(F))^{\wedge}_{\bK(Spec(k))\rightarrow \mb{F}_{\ell}} \rightarrow \mfk(Spec(\overline{F}), G)^{\wedge}_{\bK(Spec(k), G)\rightarrow \mb{F}_{\ell}}$$is an equivalence.   
 \end{Corollary}
 \begin{Proof}  The equivalence from Item \ref{fixedfield} of Proposition \ref{descentproperties} is a map of $\bK(Spec(k))$-modules. Therefore its  completion at the homomorphism $\bK(Spec(k))\rightarrow \mb{F}_{\ell}$ is an equivalence of $\bK(Spec(k))^{\wedge}_{\mb{F}_{\ell}}$-modules.   Item \ref{zeroness}  of Proposition \ref{descentproperties} shows that $\mfk(Spec(\overline{F}), G) $ satisfies the hypotheses of Lemma \ref{sphereequiv}, which gives the result.  
 \end{Proof} 

 We record another simple result. 
 \begin{Proposition} \label{transfercompletion} Let $A \rightarrow B \rightarrow C $ be a diagram of commutative ring spectra.  Suppose that $M \rightarrow N$ is a map of $B$-modules, and that $id _C \wedge _A f: C \wedge _A M \rightarrow C \wedge _AN$ is an equivalence.  Then $id_C \wedge _Bf:C \wedge_{B}  M \rightarrow  C \wedge _B N$ is an equivalence.  
 \end{Proposition}
\begin{Proof}
 Let $\mc{C}(f) $ denote the cofibre of $f$, which is a $B$-module.  The hypothesis is equivalent to the assertion that  $C \wedge _A\mc{C}(f) \simeq *$. The module $C \wedge _B\mc{C}(f)$ is also equivalent to a point, since it can be expressed as 
 $$ C \wedge _{C \wedge _A B} (C \wedge _A \mc{C}(f)) 
 $$
 This implies that $id_C \wedge_B f$ is an equivalence.  
\end{Proof}

Finally, we  reinterpret the completion $P^{\wedge}_{S^0 \rightarrow \mb{F}_{\ell}}$ in terms of inverse systems of Moore spectra.  The Moore spectrum is the cofibre of the map $S^0 \xrightarrow{\times n} S^0$, so $H_0(M_n) \cong\mb{Z}/n \mb{Z}$.  It is easy to assemble the Moore spectra $M_{\ell ^j}$ into an inverse system of spectra so that the map $H_0(M_{\ell ^i}) \rightarrow H_0(M_{\ell  ^{i-1}})$ is the projection $\mb{Z} / \ell^{i}\mb{Z} \rightarrow \ell^{i-1}\mb{Z} $, and this can be done in an essentially unique way. 

 We apply this to the $G$-equivariant  $\bK$-theory, where $G$ is the absolute Galois group of a geometric field. 

 \begin{Corollary}\label{trivialcorollary}
     Let $F$ be a geometric field of characteristic zero, with $\overline{F}$ its algebraic closure, and  with $G$ its absolute Galois group.  Then the natural map $\mfk(\overline{F},G) ^{\wedge }_{{S^0} \rightarrow \mb{F}_{\ell}} \rightarrow \mfk(\overline{F},G)^{\wedge}_{\bK \rightarrow \mb{F}_{\ell}}$ is an equivalence, where $\bK$ denotes $\bK(Spec(\overline{F}),G) $.  
 \end{Corollary}
 \begin{Proof} Item \ref{zeroness} of Proposition \ref{descentproperties} shows that the hypotheses of Lemma \ref{sphereequiv} are satisfied.  
     \end{Proof}
 \begin{Lemma}\label{Moorestuff}  Let $A$ be a commutative ring spectrum, equipped with a homomorphism $A \rightarrow \mb{F}_{\ell}$. For any $A$-module spectrum $P$, the natural map 
 $$P^{\wedge}_{\mb{F}_{\ell}} \rightarrow \underset{i}{holim} (P \wedge M_{\ell ^i})^{\wedge}_{\mb{F}_{\ell}}$$
is an equivalence of spectra, where $M_{\ell ^i}$ is the $\ell ^i$-Moore spectrum defined above.   
      \end{Lemma}
      \begin{Proof} It is easy to see that it will suffice to prove that the  map 
      $$ T_A(P; \mb{F}_{\ell}) \rightarrow \underset{i}{holim}\hspace{.03cm}T_A(P \wedge M_{\ell ^i})
      $$
      is an equivalence.  There is a cofibration sequence of pro-spectra 
      $$ \{K_i \}_i \rightarrow S^0 \xrightarrow{p} \{ M_{\ell ^i}\}_i
      $$
      where $S^0$ denotes the constant pro-spectrum with value $S^0$, and where $\{ K_i \}_i $ denotes the pro-spectrum 
       $$ \cdots \xrightarrow{\times \ell}  S^0 \xrightarrow{\times \ell}  S^0 \xrightarrow{\times \ell}  S^0
 $$
      The map in the statement of the Lemma is $(\times \ell)^{\wedge}_{\mb{F}_{\ell}}$, so to prove the result, it will suffice to prove that $\underset{i}{holim}\hspace{.05cm} T_A(P \wedge K_i) \simeq *$. This is immediate since multiplication by $\ell$ is trivial on $T_A(M;\mb{F}_{\ell})$ \end{Proof}

 \section{A Calculation} \label{Calculation}

 In this section, we will study a class of pro-continuous $G$-schemes, in the sense of Definition \ref{procontinuous},  of profinite groups $G$ on pairs $(B,Y)$, where $B$ is a $G$-reachable $A$ algebra for a a $k$-algebra $A$, and $Y \subseteq Spec(A)$ is a $G$-invariant  affine  subscheme, for which we can prove that $\mb{F}_{\ell} \wedge _{\bK(Spec(k),G)}\mfk_{ Y}(B, G)  \simeq *$.  This class of actions and the  vanishing result  will then be useful for our computations of derived completions.  We'll need to introduce some terminology to define it.  

As always, let $k$ denote an algebraically closed field.  
\begin{Definition} \label{unidefinition} For $i$ a positive integers, let  $P^i$ denote the polynomial algebra $k[t_1, \ldots , t_i]$.   Let   $\mathfrak{M}$ denote the left  filtering partially ordered set whose elements are the  positive integers,and where $m \leq  n$ if $n|m$, so that there is morphism $m \rightarrow n$. For all  positive integers $i$ and  $m$, we let 
$$P^i _m = k[t^{\frac{1}{m}}_1, \ldots ,  t_i ^{\frac{1}{m}}]$$
so $P^i = P^i _1$.  
We note that whenever $m|n$, there is an inclusion $P^i _m \hookrightarrow P^i _n$. These inclusions  define a $P^i $-algebra-valued functor $\mu$  on $\mathfrak{M}^{op}$, whose colimit we denote by $P_{\infty}^i $. It is easy to observe that $P_{\infty}^i $ is  a reachable $P^i$-algebra, and that $\mu$ takes its values in the strongly flat $P^i$ subalgebras under $P^i _{\infty}$.  We refer to $P_{\infty}^i $ as the {\em Puiseux algebra} of dimension $i$ over $k$.  

We also define the profinite group $\Gamma ^i = \Sigma _i \ltimes \hat{\mu}^i$, where $\hat{\mu}$ denotes the profinite group whose  parameter category is $\mathfrak{M}$ and whose corresponding functor is given by $n \rightarrow \mu_n$, where $\mu _n$ denotes the
$n$-th roots of unity. 
Note that  $\Gamma ^i$ acts on $P^i_m $  coordinatewise via multiplication on the coordinates for $\mu _m^i$,  by permutation of coordinates for $\Sigma _i$, and on $P^i $ simply by  permutation of coordinates.  The  $\Gamma ^i$-action on $P_{\infty}^i$ is continuous, and $P_{\infty}^i $ is a $\Gamma ^i$-reachable $P^i$-algebra. We will  consider the following  more general situation.  Let $\Gamma$ be a profinite group, $f:\Gamma \rightarrow \Gamma ^i$ a  continuous homomorphism, and $\varphi:\Gamma \rightarrow Aut_A(B)$ a continuous homomorphism, where $B/A$ is a finite faithfully flat ring  extension, the map $Spec(B) \rightarrow Spec(A) $ is surjective so the extension $A \subseteq B$ is smooth, and $A$ is Noetherian.   This data defines a continuous action of $\Gamma$ on  the algebra $P^i _{\infty} \otimes _k B$, and it is easy to verify that it is a $\Gamma$-reachable $P^i \otimes _k F$-algebra. Definition \ref{bigdefinition} defines a pro-continuous $\Gamma$-pro-scheme $Spec(P^i \otimes _k B,P^i_{\infty}\otimes _k B)$ attached to this algebra, which we will denote by $\mc{P}^i(B) $.  We are then free to evaluate   $\mfk(-,\Gamma )$ on $\mc{P}^i(E)$ to obtain a $\bK(Spec(k),\Gamma )$-module spectrum $\mfk(\mc{P}^i(B),\Gamma)$. 
 
\end{Definition}
We will now prove the following theorem. 
\begin{Theorem} \label{vanishing} With notation as in Definition \ref{unidefinition}, we have
$$ {\mb{F}}_{\ell} \wedge _{\bK(Spec(k), \Gamma)}
    \mfk_{ \{ 0 \}\times _k Spec(A)}(\mc{P}^i(B) , \Gamma) \simeq  * 
$$
    where $\{ 0 \}$ denotes the origin in $Spec(P^i) $.  $\mfk _{ \{ 0 \}\times _k Spec(A) }$ is defined in Definition \ref{pullbackschemes}. 
\end{Theorem}
\begin{Proof} 

Because $B$ is finite and faithfully flat over $A$, it is immediate that the full subcategory of $SF(P^i _{\infty} \otimes _k B)$ consisting of the algebras $P^i_{m}  \otimes _k B$ is cofinal, and hence that it suffices to prove that 
$$ \underset{m \in \mathfrak{M}^{op}}{colim} \hspace{.05cm}\mb{F}_{\ell}\wedge _{\bK(Spec(k),G)} \mfk _{\pi ^{-1}_{P^i_m \otimes B}(\{ 0 \} \times Spec(A) )}(Spec(P^i_m \otimes B),G) \simeq *
$$

For economy of notation, we let $\mc{G}(m) = \mfk _{\pi ^{-1}_{P^i_m\otimes B}(\{ 0 \} \times Spec(A)} (Spec(P^i _m \otimes _k B),G)$. We observe that each $\mc{G}(m) $ is a module over the commutative ring spectrum $\bK(Spec(B), G)$, and that $\mc{G}$ is a functor from $\mathfrak{M}^{op}$ to the category of $\bK(Spec(B), G)$-modules.   Because the partially ordered set $\mathfrak{M}^{op}$ is  filtering, it further follows that it suffices to show that the colimit of graded abelian groups 
$$  \underset{m \in \mathfrak{M}^{op}}{colim} \hspace{.05cm}\pi _* (\mb{F}_{\ell}\wedge _{\bK(Spec(k),G)} \mc{G}(m))  $$
vanishes. 

In order to prove this, it is sufficient to prove that for every $m \in \mathfrak{M}^{op}$ there is another $m^{\prime} \in \mathfrak{M}^{op}$ so  that  (a) $m \leq m^{\prime}$ and (b) the map of spectra 
$$ \mb{F}_{\ell}\wedge _{\bK(Spec(k),G)} \mc{G}(m)  \rightarrow \mb{F}_{\ell}\wedge _{\bK(Spec(k),G)} \mc{G}(m^{\prime}) 
$$
 induces the trivial map on homotopy groups.  We will show that the selection $m^{\prime} = \ell m$ satisfies these conditions.  It is trivial that it satisfies (a). We define $\alpha_m$ to be the map $\mc{G}(m \leq m \ell)$, and $\beta _m$ to be $\mb{F}_{\ell} \wedge _{\bK (Spec(k),G)}\alpha _m$.  We will show that $\pi _* \beta _m$ vanishes, which will prove the theorem.

We observe that $\bK(Spec(B),G)$ is a commutative $\bK(Spec(k),G)$-algebra  spectrum, that all the spectra $\mc{G}(m)$ are $\bK(Spec(B),G)$-modules, and that the maps $\mc{G}(m) \rightarrow \mc{G}(m^{\prime})$ are maps of $\bK(Spec(B),G)$-modules. Moreover, because of the smoothness and the fact that  $B$ is Noetherian  (since $A$ is), and using item \ref{relequiv} of Proposition \ref{kproperties}, we have natural equivalences of $\bK(Spec(B), G)$-modules 
$$ \bK(Spec(B),G) \xrightarrow{\sim}\mfk(Spec(B),G) \xrightarrow{\sim}  \mfk _{\pi ^{-1}_{P^i_n \otimes B}(\{ 0 \} \times Spec(A) )}(Spec(P^i_n \otimes B),G)= \mc{G}(n)
$$
for all $n$.  
We'll denote the equivalence $\mfk(Spec(B),G) \rightarrow \mc{G}(n) $ by $\varepsilon_n$, and  construct a homotopy inverse $p_n$ to it.  Considering the inclusion $B \rightarrow P^i_n \otimes _k B $, we note that if $M$ is a $P^i_n \otimes _k B$-module supported on $\{ 0 \} \times Spec(B)$, then the $B$-module obtained by restricting scalars along the map $B \rightarrow P^i_n \otimes _k B$ is finitely generated.  We therefore have  a map of $\bK(Spec(B),G)$-module spectra from $ \mc{G}(n) $ to $ \mfk (Spec(B),G) $, and it is readily checked that $p_n$ is a homotopy inverse to $\varepsilon_n$ for all $n$.  We will now examine the composite 
$$ \mfk(Spec(B) ,G) \xrightarrow{ \varepsilon _m} \mc{G}(m) \xrightarrow{\alpha _m} \mc{G}(m \ell)  \xrightarrow{p_{m \ell} } \mfk(Spec(B),G)
$$
Since $\varepsilon_m $ and $p_{m \ell}$ are both equivalences, it follows that the vanishing of $\pi _*\beta _m$ is equivalent to the vanishing of $\pi _* (p_{m \ell }\circ\alpha _m \circ \varepsilon_m ) $. We now describe the functors inducing $\varepsilon_m$, $\alpha _m$, and $p_{m \ell}$.  

\begin{enumerate}
    \item{{\bf $\bm{\varepsilon _m }$: } Given a finitely generated $B$-module $M$, $\varepsilon _m$ constructs  the $P^i_m \otimes _k B$-module $k \otimes _k M$, where $P^i_m \otimes _k B$ acts via $p \otimes b \cdot (\kappa \ \otimes m)  = p\, \cdot \kappa \otimes b \cdot m$, for $\kappa \in k$ and $m \in m. $}
    \item{{\bf $\bm{\alpha _m}$:} Given a $P^i_m \otimes _k B$-module $M$, $\alpha_m $ assigns to it the $P^i_{m \ell} \otimes _k B$-module $(P^i_{m \ell}  \otimes B) \otimes _{P^i _m \otimes _k B} M$.}
    \item{{\bf $\bm{p_{m \ell}}$: } Given a $P^i_{m \ell} \otimes _k B$-module $M$ supported on $\{ 0 \}\times Spec(B)$, $p_{m \ell}$ assigns  to it $M$ regarded as a $B$-module by restriction of scalars to $B \subseteq P^i_{m \ell} \otimes _k B$. The result is clearly finitely generated.   }
\end{enumerate}
We claim the composite functor $p_{m \ell} \circ \alpha _m \circ \varepsilon _m $ is given by $M \rightarrow (P^i_{m \ell} \otimes _{P^i _{m}}k) \otimes _k M$. This follows since there is a standard isomorphism 
$$    (P^i_{m \ell}\otimes _k B )\otimes _{P^i _m \otimes _k B}(k \otimes _k M) \rightarrow (P^i_{m \ell} \otimes _{P^i _m} k ) \otimes _k M 
$$
of $G$-equivariant $B$-modules.  The $G$-module $P^i_{m \ell} \otimes _{P^i_m} k $ is a $k$-linear representation $\rho $ of finite dimension $\ell ^i$, and therefore represents an element  $[\rho]\in R[G]$ with trivial mod-$\ell$ augmentation.  It now follows that the homotopy class $p_{m \ell} \circ \alpha _m \circ \varepsilon _m $ is represented by $$[ \rho ] \wedge id : \bK(Spec(k),G) \wedge _{\bK (Spec(k),G) }\mfk (Spec(B),G)  \rightarrow \bK(Spec(k),G) \wedge _{\bK (Spec(k),G) }\mfk (Spec(B),G) $$
It follows that $\mb{F}_{\ell} \wedge_{\bK(Spec(k),G)} (p_{m \ell} \circ \alpha _m \circ \varepsilon _m) $ is null-homotopic since $\mb{F}_{\ell} \wedge ([\rho] \wedge id) \simeq (\mb{F}_{\ell} \wedge [\rho])  \wedge id) $, and $\mb{F}_{\ell} \wedge [ \rho] \simeq *$ since the mod-$\ell$ augmentation of $\rho$ is zero. 
\end{Proof}

The following corollary will be useful in Section \ref{Bextensions}
\begin{Corollary} \label{subgroups}
Given the situation in Theorem \ref{vanishing}, and a subgroup $\Gamma _0 \subseteq \Gamma $ of finite index.  Then  
$$ {\mb{F}}_{\ell} \wedge _{\bK(Spec(k), \Gamma)}
    \mfk_{ \{ 0 \}\times _k Spec(A)}(\mc{P}^i(B) , \Gamma_0 ) \simeq  *
$$
    \end{Corollary}
    \begin{Proof}
        Let $i_{\Gamma _0}^{\Gamma}(A)$ (respectively $i_{\Gamma _0}^{\Gamma}(B)$)  denote the algebras of $\Gamma _0$ invariant functions from $\Gamma$ to $A$ (respectively $B$), where $\Gamma_0$ acts on $\Gamma $ on the right, and on the target $B$.  The left action of $\Gamma$ on $\Gamma$ defines  $\Gamma$-actions on $i_{\Gamma _0}^{\Gamma}(A)$ and $i_{\Gamma _0}^{\Gamma}(B)$, with an equivariant inclusion 
        $i_{\Gamma _0}^{\Gamma} (A) \hookrightarrow i_{\Gamma _0}^{\Gamma} (B)$.  Applying the theorem with $A$ and $B$ replaced by $i_{\Gamma _0}^{\Gamma}(A)$ and $i_{\Gamma _0}^{\Gamma}(B)$ gives the result.  
    \end{Proof}

 \section{$K$-theory of Puiseux actions}\label{Bextensions} 

We will assume we are in the situation of Theorem \ref{vanishing}, where we are given a profinite group $G$, a continuous homomorphism $f: G \rightarrow \Gamma^i$ and $B$ is a finite faithfully flat $A$-algebra with compatible $G$-action.   In this section, the algebra $B$ will be a field containing $k$, equipped with a continuous homomorphism  with finite image from a $G$ to the automorphism group of $E$ over $k$. We let $F$ be the fixed field of the $G$-action, so $E/F$ is a finite Galois extension. We include the possibility that $E = k$, with the group action on $k$ assumed to be trivial.  For any $G$-invariant subscheme $Y \subseteq Spec(P^i)$, we define the pro-continuous $G$-scheme $\mc{P}^i(E, Y) = Spec(P^i \otimes_k E, P^i _{\infty} \otimes _k E; Y) $ as in Definition \ref{pullbackschemes}. We  define $\mb{A} = Spec(P^i)$ and let  $\mb{T} \subseteq \mb{A} $ denote the subscheme obtained by removing the coordinate hyperplanes in $\mb{A}$. In this section we will prove the following theorem. 

\begin{Theorem}\label{mainfinal}
    The map of $\bK(Spec(k), G)$-module spectra 
    $$ \mb{F}_{\ell} \wedge _{\bK (Spec(k),G)} \mfk (\mc{P}^i(E,\mb{A}),G) \rightarrow \mb{F}_{\ell} \wedge _{\bK (Spec(k),G)} \mfk (\mc{P}^i(E,\mb{T}),G)
    $$
    is an equivalence of spectra.  It follows from item \ref{keyprop} of Proposition \ref{properties} that the map 
    $$\mfk (\mc{P}^i(E,\mb{A}),G)^{\wedge}_{\mb{F}_{\ell}} \rightarrow  \mfk (\mc{P}^i(E,\mb{T}),G)^{\wedge}_{\mb{F}_{\ell}}
    $$
    is an equivalence.  It is a straightforward consequence of the homotopy property that $\mfk(Spec(E),G) \rightarrow \mfk (\mc{P}^i(E,\mb{A}),G)$ is an equivalence, so we may conclude that the map 
    $$ \mb{F}_{\ell} \wedge _{\bK (Spec(k),G)} \mfk (Spec(E),G) \rightarrow \mb{F}_{\ell} \wedge _{\bK (Spec(k),G)} \mfk (\mc{P}^i (E,\mb{T}),G)
    $$
\end{Theorem}

\begin{Remark} \label{rem.fundfiltration}
In view of the above theorem, the object $\{\mc{P}^i _{\mc{\mb{T}}}|i\}$ should be viewed as
model of the universal bundle over a model of the classifying space for the profinite group $G$: as such one may denote
it by ${\rm EG}^C$ as discussed in \cite[section 4(ii)]{royg4}. (Observe also that each $\mc{P}^i _{\mc{\mb{T}}}$ is a smooth scheme.)

 Therefore, the last theorem may be viewed as a derived completion theorem for $\bG(\quad, G)$-theory
with respect to the action of the profinite group $G$ on the inverse system of
affine schemes $\{Spec\mc(A)_i|i\}$. As such, the main statement of the theorem could be rewritten as:
\[\mfk(\mc{A}, G)^{\wedge}_{\mb{F}_{\ell}} \cong \mfk({\rm EG}^C\times_{k} \mc{A}, G) ^{\wedge}_{\mb{F}_{\ell}}. \]
\end{Remark}
{\bf Proof:} 
We let $\mb{S} = \mb{A}- \mb{T}$ and now have a localization sequence   
\begin{equation}  \mfk _{\mb{S}} (\mc{P}^i(E,\mb{A}), G) \rightarrow \mfk (\mc{P}^i(E, \mb{A}), G )  \rightarrow \mfk (\mc{P}^i(E, \mb{T}),  G)
\end{equation} from Proposition \ref{pullbacklocalization}
so it will suffice to prove that  
\begin{equation} \label{reductionone} \mb{F}_{\ell} \wedge _{\mfk(Spec(k), G)}\mfk _{\mb{S}}(\mc{P}^i(E,\mb{A}), G ) \simeq * .\end{equation} There is a decreasing filtration of closed subschemes 
\begin{equation}\label{filtration}  F_i \mb{S}\subseteq F_{i-1} \mb{S}\subseteq \cdots \subseteq F_1\mb{S} = \mb{S} \end{equation}
where $F_j \mb{S}$ denotes the  image in $Spec(P^i)$ of the union of all non-degenerate $j$-fold intersections of coordinate hyperplanes.  We obtain the corresponding filtration of $\bK(Spec(k), G)$-module spectra 
\begin{equation}\label{firstsequence} \mfk_{F_i \mb{S}}(\mc{P}^i(E,\mb{A}), G) \subseteq \mfk_{F_{i-1} \mb{S}}(\mc{P}^i(E,\mb{A}), G )\subseteq \cdots \subseteq \mfk_{F_1 \mb{S}}(\mc{P}^i(E,\mb{A}), G) = \mfk_{\mb{S}}(\mc{P}^i(E,\mb{A}), G) \end{equation}

Moreover, the cofibers of the maps  
$$ \mfk_{F_{j+1} \mb{S}}(\mc{P}^i(E,\mb{A}), G) \rightarrow \mfk_{F_j \mb{S}}(\mc{P}^i(E,\mb{A}), G)
$$
are equivalent to 
\begin{equation} \label{beginning} \mfk_{F_j\mb{S}  -F_{j+1}\mb{S} }(\mc{P}^i(E,{\mb{A}-F_{j+1}\mb{S}}), G) \end{equation}
again by Proposition \ref{pullbacklocalization}. 

 We observe that there is a disjoint union decomposition 
\begin{equation} \label{subsets1}
F_j \mb{S} - F_{j+1}\mb{S} \cong  \coprod _{\#(S) = j}\left[ \bigcap _{s \in S} {\mathbb{H}}^s-\bigcup_{t \notin S}{\mb{H}}^t \cap (\bigcap_{s \in S}{\mb{H}}^s) \right]
\end{equation}
where $S \subseteq \mc{I} = \{ 1, \ldots , i \}$ ranges over the subsets of cardinality exactly $j$, and  $\mb{H}^s \subseteq \mb{A}$ denotes the $s$-th hyperplane in $\mb{A}$. Let $\Pi_j$ denote the collection of subsets $S \subseteq \mc{I}$ of cardinality exactly $j $, and for $S \in \Pi_j$, we let $\mb{H}(S)$ denote $\bigcap _{s \in S} {\mathbb{H}}^s-\bigcup_{t \notin S}{\mb{H}}^t \cap (\bigcap_{s \in S}{\mb{H}}^s) $.  The subsets $\mb{H}(S)$ are both open and closed,  as well as disjoint in $F_j \mb{S} - F_{j+1}\mb{S}$. We also observe  that  if we coordinatize  $\mb{A}^i(k)$ as $\mb{A} ^j(k)  \times _k \mb{A} ^{i-j}(k ) $, where $\mb{A}^i(k) $ denotes the $i$-dimensional Puiseux scheme on the $i$ variables belonging to $S$, and $\mb{A}^{i-j}(k) $ corresponding to  $\mc{I}-S$, then $\mb{H}(S) = \{ 0_S \}  \times \mb{T}^{i-j}$.  
 To analyze the equivariant $\mfk$-theory of these schemes, we introduce the following lemmas.  

\begin{Lemma}\label{disjoint} Let $X$ be a scheme  acted on by a finite group $G$, and $Y$ a closed subscheme invariant under the $G$-action.  We suppose that $Y$ decomposes into a disjoint union of  connected components $Y = \coprod _{\alpha \in A} Y_{\alpha}$, so  $G$ acts on the set of components of $Y$ and correspondingly on the set $A$. Let $B$ denote the set of orbits of the $G$-action on $A$, and for $\beta \in B$, we denote by $Y(\beta)$ the disjoint union $\coprod _{\alpha \in \beta} Y_{\alpha}$. This decomposition is equivariant.    

Then we have 
$$  \mfk_{Y}(X. G) \cong  \prod_{\beta \in B} \mfk_{Y(\beta)}(X, G)
$$
    
\end{Lemma}
\begin{Proof}
We consider the following diagram. 
\begin{equation}
    \xymatrix{{\prod _{\beta \in B}\mfk _{Y(\beta)}(X, G)}  \ar@<1ex>[r] \ar@<1ex>[d]^{\simeq} &{\mfk_{Y}(X, G)} \ar@<1ex>[d]^{\simeq} \\
    {\prod _{\beta \in B}\mfk(Y(\beta), G)} \ar@<1ex>[r]^{\simeq} &   {\mfk(Y, G)} .}
\end{equation}

    The upper horizontal arrow is  a product of maps induced by inclusions of exact subcategories of more  restrictive support to a subcategory of less restrictive support.  
    The vertical arrows are equivalences arising from  from devissage arguments. The lower horizontal arrow is an equivalence since $Y_0 = \coprod _{\beta \in B} Y_0(\beta)$. 
\end{Proof}
\begin{Lemma} \label{disequivariant}
Let $Y = \coprod_{\alpha \in A} Y_{\alpha}$ be a scheme, and let $G$ be a finite group acting on $Y$, permuting the subschemes $Y_{\alpha}$ transitively.  For each $\alpha$, we suppose we are given a closed subscheme  $Z_{\alpha } \in Y_{\alpha}$, and suppose that the set $Z = \coprod _{\alpha}Z_{\alpha}$ is invariant under the $G$-action.  Let $\alpha \in A$, and let $H_{\alpha}$ denote the stabilizer of the element $\alpha$. Then 
$$\mfk(Y, G) \cong \mfk(Y_{\alpha}, H_{\alpha})
$$
and 
$$\mfk_{Z}(Y, G) \cong \mfk_{Z _{\alpha}}(Y_{\alpha}, H_{\alpha})
$$
for all $\alpha$. 
\end{Lemma}
\begin{Proof} This is straightforward using \cite{brion} and \cite{thomason}. 
\end{Proof} 
\begin{Lemma}\label{restriction}
    Let $X$ be a scheme acted on by a finite group, and let $Y$ be a closed subscheme.  Suppose that $X \hookrightarrow \overline{X}$ is an embedding, so that $Y$ is also closed in $\overline{X}$. Then the map $\mfk_Y(\overline{X},G) \rightarrow \mfk_Y(X,G)$ is an equivalence.  
\end{Lemma}
\begin{Proof}
    Devissage shows that spectra are equivalent to $\mfk^G(Y, G)$. 
\end{Proof}

We first apply Lemmas \ref{disjoint} and \ref{disequivariant} to the following situation.  The ambient scheme $X$ from Lemma \ref{disjoint} will be $\mb{A} - F_{j+1}\mb{S}$, and $Y$ will be $F_j \mb{S}- F_{j+1} \mb{S}$.  The subschemes $\mb{H}(S)$  cover $F_j \mb{S}- F_{j+1} \mb{S}$ disjointly and are therefore open and closed.   Moreover, the $G$-action permutes them via the action of $G$  on the collection of subsets of $\mathcal{I}$ of cardinality exactly  $j$. As noted above, the scheme $\mb{H}(S) $ is also affine, and decomposes as $\{ 0_S \} \times \mb{T}^{i-j}$ when $\mb{A}^i(k) $ is decomposed as $\mb{A}^j(k) \times \mb{A}^{i-j}(k)$, with $\mb{A}^j(k) $ the affine space on the variables in $S$ and $\mb{A}^{i-j}(k) $ on the variables  in $\mc{I}-S$ 
The action of the stabilizer of $S$ under the $G$-action preserves the  product decomposition $\mb{A}^i(k) \cong \mb{A}^j (j) \times \mb{A}^{i-j}(k)$.   Note that the stabilizer of $S$ is in fact a subgroup of $\Sigma _j\times  \Sigma _{i-j } \subseteq \Sigma _i $, with its natural action on $\mb{A}^{j} \times \mb{A}^{i-j}$.  The indexing set $A$  in Lemma \ref{disjoint} is the collection of subsets $S \subseteq \mc{I}$ of cardinality exactly $j$, and for such a subset $S$ the corresponding disjoint union factor $Y_S$ is  $ \mb{H}(S)$.  Let $\mc{O}_j$ denote the collection of $G$-orbits of the action of $G$ on the sets of cardinality exactly $j$. For each orbit $\mathfrak{o} \in \mc{O}_j$, let $$\mb{H}_{\mathfrak{o}} = \bigcup _{S \in \mathfrak{o}} \mb{H}(S)$$  Then the   conclusion of Lemma \ref{disjoint} is that we have a decomposition
\begin{equation} \label{initdecomp} \mfk_{F_{j}\mb{S}- F_{j+1} \mb{S}}(\mb{A}- F_{j+1}\mb{S}, G) \cong \prod _{\mathfrak{o} \in \mc{O}_j}\mfk_{\mb{H}_{\mathfrak{o}}}(\mb{A}- F_{j+1} \mb{S}, G)
\end{equation}

Next, for each orbit $\mathfrak{o} \in \mc{O}_j$, let $S_{\mathfrak{o}}$ denote an orbit representative, and let $G_{\mathfrak{o}}$ denote the stabilizer of $S_{\mathfrak{o}}$.  Then Lemma \ref{disequivariant} asserts that 
\begin{equation} \label{seconddecomp} \mfk_{\mb{H}_{\mathfrak{o}}}(\mb{A}- F_{j+1} \mb{S}, G) \cong \mfk_{\mb{H}(S_\mathfrak{o})}(\mb{A}- F_{j+1} \mb{S}, G_{\mathfrak{o}})
\end{equation}

Next, let $\overline{S}_{\mathfrak{o}}$ denote $\mathcal{I}- S_{\mathfrak{o}}$.  For any $s \in \overline{S}_{\mathfrak{o}}$, let $V(s)$ denote the subscheme of $A$ defined by $t_s = 0$. 
We note that $\{ 0 _S\} \times \underset{s \in \overline{S}_{\mathfrak{o}}}{\cup}V(s) \subseteq  F_{j+1}\mb{S} $, viewed in a coordinate system  respecting the decomposition $\mb{A}^i\cong \mb{A}^{j}(k) \times \mb{A}^{i-j}(k)$. It follows by Lemma \ref{restriction} that the map
\begin{equation}\label{fundequiv} \mfk_{\mb{H}(S_\mathfrak{o})}(\mb{A}- F_{j+1} \mb{S}, G_{\mathfrak{o}}) \rightarrow \mfk_{\{ 0 _S\} \times \mb{T}^{i-j}}(\mb{A}^{j}\times \mb{T}^{i-j}, {G}_{\mathfrak{o}})
\end{equation}
induced by that inclusion is an equivalence.  
This map  is an equivalence of $\mfk(Spec(k), G)$-modules. 

The discussion around (\ref{firstsequence}) and (\ref{beginning}) above showed that the cofibers of the map of spectra 
$$ \mfk_{F_{j+1}\mb{S}}(\mc{P}^i(E, \mb{A}), G) \rightarrow \mfk_{F_{j}\mb{S}}(\mc{P}^i(E, \mb{A}), G) 
$$
are equivalent to the spectra $\mfk_{F_j\mb{S} - F_{j+1}\mb{S}}(\mc{P}^i (E
,{\mb{A}-F_{j+1} \mb{S}}), G)$.

The constructions in (\ref{initdecomp}), (\ref{seconddecomp}), and (\ref{fundequiv}) above now translate  to corresponding constructions for $\mc{P}^i $, namely 
\begin{equation} \mfk_{F_{j}\mb{S}- F_{j+1} \mb{S}}(\mc{P}^i(E ,{\mb{A}- F_{j+1}\mb{S}}), G) \cong \prod _{\mathfrak{o} \in \mc{O}_j}\mfk_{\mb{H}_{\mathfrak{o}}}(\mc{P}^i(E,{\mb{A}- F_{j+1} \mb{S}}), G),
\end{equation}
\begin{equation}
    \mfk_{\mb{H}_{\mathfrak{o}}}(\mc{P}^i(E,{\mb{A}- F_{j+1} \mb{S}}), G) \cong \mfk_{\mb{H}(S_\mathfrak{o})}(\mc{P}^i(E,{\mb{A}- F_{j+1} \mb{S}}), G_{\mathfrak{o}}),
\end{equation}
and 
\begin{equation}  \mfk_{\mb{H}(S_\mathfrak{o})}(\mc{P}^i(E,{\mb{A}- F_{j+1} \mb{S}}), G_{\mathfrak{o}}) \rightarrow \mfk_{\{ 0 _S\} \times \mb{T}^{i-j}}(\mc{P}^i(E,{\mb{A}^{j}\times \mb{T}^{i-j}}), {G}_{\mathfrak{o}}).
\end{equation}
It follows that to prove the theorem, it suffices to prove that 
$$ \mb{F}_{\ell} \wedge _{\mfk(Spec(k), G)}\mfk_{\{ 0 _S\} \times \mb{T}^{i-j}}(\mc{P}^i(E,{\mb{A}^{j}\times \mb{T}^{i-j}}), {G}_{\mathfrak{o}}) \cong *
$$

The parameter category for the Ind-spectrum $\mfk_{\{ 0 _S\} \times \mb{T}^{i-j}}(\mc{P}^i(E, {\mb{A}^{j}\times 
\mb{T}^{i-j}}), {G}_{\mathfrak{o}}) $ admits a final  subcategory equivalent to the  by the category $\mathfrak{M}^{op}$, with the restriction of the spectrum valued functor  given by $ m \rightarrow \mfk_{\{ 0 _S\} \times \mb{T}^{i-j}_m}(\mc{P}^i(E, {\mb{A}^{j}_m\times \mb{T}^{i-j}_m}), {G}_{\mathfrak{o}})$, where $\mb{A}^j_m$ and $\mb{T}^{i-j}_m $ are defined by 
$$ \mb{A}^j_m = p_m^{-1}(\mb{A}^j ) \mbox{ and } \mb{T}^{i-j}_m = p_m^{-1}(\mb{T}^{i-j} )
$$
where $p_m: \mb{A}^i_m \rightarrow \mb{A}^i$
is the projection.  

We may construct a two parameter prospectrum with parameter category $\mathfrak{M}^{op} \times \mathfrak{M}^{op}$, with defining functor
$$ (m,m^{\prime}) \rightarrow \mb{F}_{\ell} \wedge _{\bK(Spec(k),G)}\mfk_{\{ 0 _S\} \times \mb{T}^{i-j}_{m^{\prime}}}(\mc{P}^i(E, {\mb{A}^{j}_m\times \mb{T}^{i-j}_{m^{\prime}}}), {G}_{\mathfrak{o}})
$$
It is immediate that the colimit spectrum over $\mathfrak{M}^{op} \times \mathfrak{M}^{op}$ is equivalent to the one-variable version, via the cofinality of the diagonal inclusion $\mathfrak{M}^{op} \rightarrow \mathfrak{M}^{op}  \times \mathfrak{M}^{op}  $. The two variable version can, however, be computed iteratively as 
$$\underset{m^{\prime}}{colim}\Big[\underset{m}{colim}\hspace{.05cm} \mb{F}_{\ell} \wedge _{\bK(Spec(k), G)} \mfk_{\{ 0 _S\} \times 
\mb{T}^{i-j}_{m^{\prime}}}(\mc{P}^i(E, {\mb{A}^{j}_m\times \mb{T}^{i-j}_{m^{\prime}}}), {G}_{\mathfrak{o}})\Big]
$$
It remains to analyze 
\begin{equation} \label{finalequiv}
\underset{m}{colim} \mbox{\hspace{.07cm}}
\mb{F}_{\ell}\wedge _{\bK(Spec(k),G)}\mfk_{\{ 0 _S\} \times 
\mb{T}^{i-j}_{m^{\prime}}}(\mc{P}^i(E, {\mb{A}^{j}_m\times \mb{T}^{i-j}_{m^{\prime}}}),  {G}_{\mathfrak{o}})
\end{equation}

with $m^{\prime}$ fixed. The affine coordinate ring of $\mb{T}_{m^{\prime}}^{i-j}$ satisfies the requirements of Theorem \ref{vanishing}, so we can conclude that the colimit (\ref{finalequiv}) is contractible.   Note that we have applied Corollary \ref{subgroups} to apply Theorem \ref{vanishing}.   \qed

 \section{Totally Torsion Free Groups}
We summarize the results of  \cite{royg1}.  A  profinite group $G$  is said to be {\em virtually finitely generated  torsion free abelian, or is a  VTF-group} if it contains a closed subgroup of finite index which is finitely generated torsion free abelian. It is shown in \cite{royg1}   that the VTF condition is equivalent to the condition that $G$ embeds as a closed subgroup of $\Gamma ^i  = \Sigma _i \ltimes \hat{\mb{Z}}^i$ for some $i$. We say a  profinite group $G$ is {\em VTF-approximable } if for every closed normal finite index subgroup $K\subseteq G$, there is a closed normal (in $G$)  subgroup $H \subseteq K$  so that $G/H$ is VTF.
We also let $\mathfrak{B}$ denote the class of torsion free finitely generated  VTF profinite groups, and say that a profinite group $G$ is $\mathfrak{B}$-approximable  if in the above definition, the quotient $G/H$ may be taken to lie in $\mathfrak{B}$. 

We say a  profinite group $G$ is {\em totally torsion free} if every closed subgroup in $G$ has torsion free abelianization.   Theorem 3.1 of \cite{royg1} is the following. 

\begin{Theorem}
    Any totally torsion free profinite group $G$ is $\mathfrak{B}$-approximable.  In particular, the absolute Galois groups of geometric fields are $\mathfrak{B}$-approximable.
\end{Theorem}

\begin{Remark} It is proved in \cite{royg1} that  separable    Galois groups for fields containing all roots of unity are totally torsion free. 
\end{Remark}
Let $\mc{X}(G) $ denote the  set of all Puiseux representations of $G$.  We define  $\Pi^G$ to be  the category  whose objects  are finite 
subsets $X = \{\rho _1, \ldots , \rho _n \} $ of $\mc{X}(G)$, and where there is a unique morphism from $X$ to $Y$ if and only if  $Y \subseteq 
X$. We define a functor $\Phi$ from $\Pi^G$ to the category of reachable pairs $(A,B) $ of $k$-algebras (see Definition \ref{specialty}) as follows, adopting the notation of Definition \ref{unidefinition}.  For any Puiseux 
representation 
$\rho: G \rightarrow \Sigma _i \ltimes \hat{\mb{Z}}^i  \subseteq Aut(P^i_{\infty})$
, we set $A(\rho) = 
 P^i (\frac{1}{t_1 \cdots t_i })$, and $B(\rho) =  P^i_{\infty}(\frac{1}{t_1 \cdots t_i })$. Both $A(\rho)$ and $B(\rho)$ are equipped with $G$-actions, and the action on $A(\rho)$ factors through a quotient $G/H$ where $H$ is a closed normal subgroups of finite index.  From Section \ref{Calculation}, it follows that $(A(\rho), B(\rho))$ is a $G$-reachable pair. The category $\Pi^G$  is a cocomplete, and has the empty set 
 as an initial object. The category of $k$-algebras with continuous $G$-action, trivial on $k$, is also cocomplete, with the colimit of a diagram $A \leftarrow B \rightarrow C$ isomorphic to $A \otimes _B C$. The functor $\Phi$ is specified up to isomorphism by the requirements that (a) $\Phi (\emptyset) = k $ and (b) $\Phi$ preserves colimits.  
 Explicitly,  $\Phi$ is now defined on objects by 
$$ \Phi (\{\rho _1, \ldots , \rho _n \} ) = (A(\{ \rho _1, \ldots , \rho _n \} ), B(\{ \rho _1, \ldots , \rho _n \} ))= (\bigotimes _{i=1}^n A(\rho _i), \bigotimes _{i=1}^nB(\rho _i)) 
$$ and is readily checked to be a functor to the category of $G$-reachable pairs.  We compose $\Phi$ with the contravariant functor $Spec(-) $ on the category of $G$-reachable pairs to obtain a functor from $\Pi^G $ to the category of pro-continuous $G$-schemes.  Since $\Pi^G$ is cofiltering, we can flatten this construction to obtain a pro-continuous $G$-spectrum, which we denote by $\mc{E}(G)$. 

 Of course $\mc{E}(G)$ is equipped with a natural map $\eta: {\mc{E}(G) }\rightarrow Spec(k)$, and therefore by pullback  a map $\eta _{\overline{F}}:  {\mc{E}(G) }\times _k Spec(\overline{F}) \rightarrow  Spec(\overline{F})  $, where $Spec(\overline{F})$ denotes pro-scheme of the spectra of  all finite extensions of $F$.     
\begin{Proposition}\label{assemblecolim}
Let $G$ denote the absolute Galois group of a geometric field of characteristic zero, containing an algebraically closed subfield $k$, and let $\overline F$ denote the algebraic closure of $F$.  
    
    Then the maps
$$ \mb{F}_{\ell} \wedge _{\bK(Spec(k),G)}\mfk(\eta ): \mb{F}_{\ell} \wedge _{\bK(Spec(k),G)}\mfk(Spec(k), G) \rightarrow  \mb{F}_{\ell} \wedge _{\bK(Spec(k),G)}\mfk({\mc{E}(G) } , G) 
   $$
   and 
   $$ \mb{F}_{\ell} \wedge _{\bK(Spec(k),G)}\mfk(\eta _{\overline{F}}): \mb{F}_{\ell} \wedge _{\bK(Spec(k),G)}\mfk(Spec(\overline{F} ), G) \rightarrow  \mb{F}_{\ell} \wedge _{\bK(Spec(k),G)}\mfk({\mc{E}(G) } \times _k Spec(\overline{F}), G) 
   $$
  are equivalences  of spectra. 
\end{Proposition}
\begin{Proof} 
We prove this for the case of $\eta _{\overline{F}}$. The case of $\eta$ is identical but is slightly simpler.  For a finite extension $E \subseteq \overline{F}$ of $F$, we can construct an analogous map of pro-continuous $G$-schemes $\eta _{E} : {\mc{E}(G)} \times _k Spec(E) \rightarrow Spec(E) $.  To prove that $\mb{F}_{\ell} \wedge _{\bK(Spec(k),G)}\mfk(\eta _{\overline{F}})$ is an equivalence, it clearly suffices to prove that $\mb{F}_{\ell} \wedge _{\bK(Spec(k),G)}\mfk(\eta _{E})$ is an equivalence for all finite extensions $E$ of $F$, since $\mb{F}_{\ell} \wedge _{\bK(Spec(k),G)}\mfk(\eta _{\overline{F}}) = \underset{E}{colim}\hspace{.05cm}\mb{F}_{\ell} \wedge _{\bK(Spec(k),G)}\mfk(\eta _{E})$. For fixed $E$, it will now suffice to prove that for any object $X$ of $\Pi ^G$, the map 
$$ \mb{F}_{\ell} \wedge _{\bK(Spec(k),G)}\mfk(Spec(E), G) \rightarrow  \mb{F}_{\ell} \wedge _{\bK(Spec(k),G)} \mfk(Spec(\Phi (X) \otimes _k E), G) 
$$
is an equivalence.  But this is the exact statement of Theorem \ref{mainfinal}, since $\Phi (X)$ is a Puiseux representation.   
\end{Proof} 

We also have the following. 

\begin{Proposition} \label{puiseuxgalois}
Let $G$ be a totally torsion free profinite group. Then the action of $G$ on the geometric points of  the pro-continuous $G$-scheme $\mc{E}(G)$ is free.  
\end{Proposition}
\begin{Proof}
    It is proved in \cite{royg1} that for any for any Puiseux representation $G$ on $\mb{P}^i$ for which corresponding the image of $G$ in $\Sigma ^i \ltimes \hat{\mb{Z}}^i$ is torsion free, the ring extension $(\mb{P}_{\infty}^i) ^G \subseteq \mb{P}^i_{\infty} $ is a Galois ring extension.  Given an object $\{ \rho_1, \ldots , \rho _n\}$ of $\Pi^G$, we will show that there is a Puiseux representation $\rho$ for which  $\Phi (\{ \rho, \rho _1, \ldots , \rho _n\})$ is a Galois ring extension.  This will prove the cofinality of the subcategory  $\Pi ^G_{Gal} \subseteq \Pi^G$ consisting of objects $x \in \Pi^G$ for which $\Phi (x)$ is a Galois ring extension.  This means that the extension 
    $( \underset{\rightarrow}{colim}\hspace{.05cm} A\{ \rho _1, \ldots , \rho _n \}, \underset{\rightarrow}{colim}\hspace{.05cm} B(\{ \rho _1, \ldots , \rho _n \})$ is Galois, which gives the result.  To construct $\rho$, we only need to show that for any arbitrary Puiseux representation $\rho _0$, there is a Puiseux representation $\rho$ so that the kernel of $\rho _0$ contains the kernel of $\rho$, and $im(\rho )$ is torsion free.  It follows from the definitions that for any Puiseux representation $\rho_0$, there is a finite index normal subgroup $N \subseteq G$ so that the restriction of $\rho _0$ to $N$ has torsion free image.  Since $N$ is finite index and normal, there is a closed normal subgroup $\overline{N}\subseteq N$ of $G$ so that $G/\overline{N}$ is VTF, and hence there is a Puiseux representation $\rho _1$ of $G/N_0$ so that the image of $\rho _1 $ is torsion free.  The Puiseux representation $\rho _0 \times \rho _1$ is now a suitable choice for $\rho$, since the kernel of $\rho$ is $Ker(\rho _0) \cap Ker (\rho _1)$ and it is easy to check that the image of $\rho$ is torsion free.  
\end{Proof}

\section{Rigidity and Proof of the Main Theorems}

The final step we will need is a rigidity result, namely that the projection $\pi _F: \mc{E}(G) \times _k Spec(\overline{F}) \rightarrow \mc{E}(G) $ induces an equivalence 
 $\mb{F}_{\ell} \wedge _{\bK(Spec(k),G)} \mfk(\pi _F, G) $. This will follow from the following statements, which are proved in \cite{royg4}. Let $G$ be a profinite group, and throughout the remainder of the paper we will let $\bK$ denote the commutative ring spectrum $\bK(Spec(k),G)$.    

 \begin{enumerate} 
 \item{There is an extension of the notion of {\em geometric classifying space $B_{geom}G$} and {\em geometric universal space $E_{geom} G$},  defined in \cite{motivichomotopy},  to profinite groups.  Of course, $E_{geom}G$ and $B_{geom}G$ are pro-schemes over a field $k$, with free $G$-actions. }
 \item{\label{hinvariance} If $X$ is a smooth (ie. every scheme in the inverse system is smooth)  pro-continuous $G$-scheme, and the action of $G$  is free on the geometric points of $X$, then the natural map of spectra 
 $\mfk (X,G) \rightarrow \mfk (E_{geom}G \times _k  X,G) $ is an equivalence of $\bK$-module spectra.  This result is  Theorem 2.4.3 (ii) of \cite{royg4}. }
 \item{\label{h2invariance}For $G = G_F$, the absolute Galois group of a geometric field $F$, and $X$ a smooth pro-continuous $G$-scheme, the natural map 
 $$  \mb{F}_{\ell} \wedge \mfk( EG \times_k  X, G) \rightarrow \mb{F}_{\ell} \wedge \mfk(EG \times_k  Spec(\overline{F}) \times_k  X) 
    $$} is an equivalence.  Here $\mb{F}_{\ell} $ denotes the mod-$\ell$ Eilenberg-MacLane spectrum.  This is Theorem 1.1.2 of \cite{royg4}. 
 
 \end{enumerate} 
 From this point on, we will set $G = G_F$, the absolute Galois group of a field $F$ of characteristic zero containing an algebraically closed subfield $k$.  

 We now have the following.   
 \begin{Proposition} \label{preliminary}
       Let $\overline{F}$ denote the algebraic closure of $F$, and let 
     $$\eta: Spec(\overline{F}) \times _k {\mc{E}(G)} 
     \rightarrow {\mc{E}(G)} $$ denote the evident projection  along the map of pro-continuous $G$-schemes $Spec(\overline{F}) \rightarrow Spec(k)$.  Then the map
    \begin{equation}\label{nearfinal} \mb{F}_{\ell} \wedge \mfk(\eta,G) : \mb{F}_{\ell} \wedge \mfk(\mc{E}(G) ,G) \rightarrow \mb{F}_{\ell} \wedge \mfk(Spec   (\overline{F}) \times _k \mc{E}(G) ,G)
    \end{equation}
     is an equivalence.  
 \end{Proposition}
 \begin{Proof}
     We have the following commutative diagram of spectra.  
     \begin{center} \vspace{-1cm}
\begin{equation} 
\begin{tikzcd}[row sep=2.5em, column sep=4em]
\mb{F}_{\ell} \wedge \mfk({\mc{E}(G)}, G) 
\arrow[r, " "] \arrow[d]
& \mb{F}_{\ell} \wedge \mfk({\mc{E}(G)}\times _k EG, G) \arrow[d,"\sim"]
 \\
\mb{F}_{\ell} \wedge \mfk(Spec   (\overline{F}) \times _k \mc{E}(G) ,G)
\arrow[r, ""] 
&  \mb{F}_{\ell} \wedge \mfk(Spec   (\overline{F}) \times _k {\mc{E}(G)}\times _k EG,G)
\end{tikzcd}
\end{equation}
\end{center}
The horizontal arrows are induced by the  projection $\mc{E}(G) \times EG \rightarrow \mc{E}(G)$, and the vertical arrows by the projection $Spec(\overline{F}) \rightarrow Spec (k) $. The left hand vertical arrow is exactly the map (\ref{nearfinal}), and the right hand arrow is an equivalence because of item \ref{h2invariance} above. The horizontal arrows are equivalences due to the freeness of ${\mc{E}(G)}$ and item \ref{hinvariance} above.  The result now follows directly. 
\end{Proof}
 
 We note that these statements involve only smashing with $\mb{F}_{\ell}$ as an ordinary spectrum rather than as a smash product over the ring spectrum $\bK$.  We wish to promote the equivalence (\ref{nearfinal}) to an equivalence of smash products over $\bK$. It is easy to see that the functor $\mb{F}_{\ell} \wedge(-)$ can be reinterpreted as the smash product $\mb{F}_{\ell} \wedge _{S^0}(-)$, since a spectrum exactly the same thing as a module spectrum over the ring spectrum $S^0$. 

 \begin{Proposition} \label{finalrigidity}
     The map 
     \begin{equation}\label{nearfinal1} \mb{F}_{\ell} \wedge _{\bK}\mfk(\eta,G) : \mb{F}_{\ell} \wedge _{\bK}\mfk({\mc{E}(G)},G_F) \rightarrow \mb{F}_{\ell} \wedge _{\bK}\mfk(Spec   (\overline{F}) \times _k {\mc{E}(G) },G)
    \end{equation}
    is an equivalence of spectra.  
 \end{Proposition}
 \begin{Proof}
 We have that the map 
 \begin{equation} \mb{F}_{\ell} \wedge _{S^0}\mfk(\eta,G) : \mb{F}_{\ell} \wedge _{S^0}\mfk((\mc{E}(G),G)
 \rightarrow \mb{F}_{\ell} \wedge _{S^0}\mfk(Spec   (\overline{F}) \times _k {\mc{E}(G)},G)
    \end{equation}
    is an equivalence by Proposition \ref{preliminary}. 
     Noting that  $\bK$ is an $S^0$-algebra, we apply Proposition \ref{transfercompletion} to obtain the result.  
 \end{Proof}

   We are now able to prove the main results. 
\begin{Theorem}\label{finaltheorem}  Let $F$ be a field of characteristic zero which contains an algebraically closed subfield $k$. Let $G$ denote the absolute Galois group of $F$.  Let $\bK$ denote the commutative ring spectrum $\bK (Spec(k),G) $.     The the following three  results hold. 
\begin{enumerate}[label=(\arabic*), font=\upshape]

\item{\label{mainequiv} The natural map 
$$ \alpha: \mfk(Spec(k), G)^{\wedge}_{\bK \rightarrow \mb{F}_{\ell} } \rightarrow \mfk(Spec(\overline{F}), G)^{\wedge}_{\bK \rightarrow \mb{F}_{\ell}}
$$ 
is an equivalence.  }
\item{The groups $\pi _*\mfk(Spec(\overline{F}), G)^{\wedge}_{\bK \rightarrow \mb{F}_{\ell}}$ are naturally isomorphic to the homotopy groups of the $l$-adic completion of the spectrum $\mfk(F)$, i.e. to $K_*(F)^{\wedge}_{\ell}$. It follows from item  \ref{mainequiv} above that the graded group  $\pi_*\mfk(Spec(k), G)^{\wedge}_{\bK \rightarrow \mb{F}_{\ell}}$ is  also isomorphic to  $K_*(F)^{\wedge}_{\ell}$.   }
\item{Let $R[G]$ denote the $k$-representation ring of $G$, and let $K_*^{M}(F) ^{\wedge}_{\ell}$ denote the $l$-adically completed Milnor $K$-theory of the field $F$.  Then the homotopy groups of the derived completion of $R[G]$  at the homomorphism $R[G]\rightarrow \mb{F}_{\ell}$, (where both $R[G]$ and $\mb{F}_{\ell}$ are actual rings, not general ring spectra), are naturally isomorphic to $K_*^{M}(F) ^{\wedge}_{\ell}$.   }
\end{enumerate}
    
\end{Theorem}
\begin{Proof}
   For Item \ref{mainequiv}, consider the diagram
\begin{center} \vspace{-1cm}
\begin{equation} 
\begin{tikzcd}[row sep=2.5em, column sep=8em]
\mb{F}_{\ell} \wedge  _{\bK} \mfk(Spec(k), G) 
\arrow[r, "\sim"] \arrow[d]
& \mb{F}_{\ell} \wedge  _{\bK} \mfk({\mc{E}(G) }, G) \arrow[d, "\sim"]
 \\
\mb{F}_{\ell} \wedge  _{\bK } \mfk(Spec(\overline{F}), G)
\arrow[r, "\sim"] 
&  \mb{F}_{\ell} \wedge  _{\bK} \mfk({Spec(\overline{F}) \times _k \mc{E}(G)}, G)
\end{tikzcd}
\end{equation}
\end{center}
The right hand vertical arrow is an equivalence by Proposition \ref{finalrigidity}. The two horizontal arrows are equivalences by Proposition \ref{assemblecolim}.  It follows immediately that the left hand vertical arrow is an equivalence.

   The second statement  is a consequence of Proposition \ref{descentproperties}  and Lemma \ref{sphereequiv}. Specifically, Item \ref{zeroness} of Proposition \ref{descentproperties} shows that hypotheses of Lemma \ref{sphereequiv} are satisfied with $A = S^0$, $B = \bK$, and $M = \mfk(Spec(\overline{F}), G)$. 
   
 The third statement is verified as follows.  We let $C$ denote the group $\mb{Z}/{\ell}^{\infty} \mb{Z}$.  The inclusion $\iota: C \rightarrow k^*$  denotes the inclusion on the ${\ell}$-power roots of unity.  $\iota $ induces a map of spectra $\Sigma ^{\infty}BC_* \rightarrow \bK(Spec(k))$, and  the inclusion of $C$ into the circle group $S^1$ induces a map $\Sigma ^{\infty}BC_+ \rightarrow BS^1_+$ which is an equivalence after $\ell$-adic completion. The group $\pi _2 (BS^1_+)$, and therefore the group $\pi _2 ((BS^1)^{\wedge}_{\ell})$, admits a canonical generator $\theta$.  The element  $\iota(\theta)$  determines  an element in $\pi _2 \bK(Spec(k))^{\wedge}_{\ell}$, which we denote by $\beta$. The evident ring spectrum map $\bK(Spec(k)) \rightarrow \bK(Spec(k), G)= \bK$ gives a corresponding element in $\pi _2 (\bK)$, which we also denote by $\beta$. For any $\bK$-module $M$, we obtain a map $\Sigma ^2 M \rightarrow M$, which we also  denote by $\times \beta$.
 
 To make the argument clearer, we will denote the spectra $\mfk(Spec(k), G) $ and $\mfk(Spec(\overline{F}), G)$ by $\mfk(k)$ and $\mfk(\overline{F})$, respectively.  We will let $\mb{H}$ denote the zero-dimensional Eilenberg-MacLane spectrum for $R[G]$.  Finally, we will adopt the notations $\mfk(k; \ell ^i)$, $\mfk(\overline{F} ; \ell^i)$, and $\mb{H}(\ell ^i) $ for the result of smashing the spectra $\mfk(k), \mfk(\overline{F})$, and $\mb{H}$ with the Moore spectrum $M_{\ell ^i}$ defined in Section \ref{derivedsection}. The spectra $\mfk(k;\ell ^i), \mfk(\overline{F}; \ell ^i)$, and $\mb{H} (\ell ^i) $ each assemble into pro-spectra parametrized by the exponent $i$.  They are all $\bK$-module spectra, and the maps in the pro-objects are maps of $\bK$-modules. The notation $(-)^{\wedge}_{\mb{F}_{\ell}}$ will denote the functor $(-)^{\wedge}_{\bK \rightarrow \mb{F}_{\ell}}$. 

 Our goal is to analyze the diagrams 
 \begin{center} \vspace{-1cm}
\begin{equation} \label{definediagram}
\begin{tikzcd}[row sep=2.5em, column sep=6em]
\Sigma ^2 \mfk(k; \ell ^i)^{\wedge}_{\mb{F}_{\ell}} \arrow[d,"\sim"]
\arrow[r, "\times \beta"] 
& \mfk(k; \ell ^i)^{\wedge}_{\mb{F}_{\ell}} \arrow[d,"\sim"] \arrow[r] \arrow[r] & \mc{C}(k; \ell ^i ) \arrow[d]
 \\
\Sigma ^2 \mfk(\overline{F} ; \ell ^i)^{\wedge}_{\mb{F}_{\ell}}
\arrow[r, "\times \beta"] 
&   \mfk(\overline{F}; \ell ^i)^{\wedge}_{\mb{F}_{\ell}} \arrow[r] &\mc{C}(\overline{F};\ell ^i )
\end{tikzcd}
\end{equation}
\end{center}
where the spectra $\mc{C}(k; \ell ^i)$ and $\mc{C}(\overline{F}; \ell ^i)$ are defined by the requirements that the rows in diagram (\ref{definediagram}) are cofibration sequences. We note that the diagrams (\ref{definediagram}) for various values of $i$ assemble into a diagram of pro-spectra.   The two left vertical arrows are equivalences by item \ref{mainequiv} of the theorem, and it follows that the right hand vertical arrow is also an equivalence.  It remains to identify the spectra $\mc{C}(k ; \ell^i)$ and $\mc{C}(\overline{F} ; \ell^i)$. 

We first analyze the cofibre of $\times \beta$ on $\mfk(k;\ell ^i)^{\wedge}_{\mb{F}_{\ell}}$. We have  weak equivalences 
\begin{equation} \label{bottequiv} \mfk(k; \ell ^i )\overset{\sim}{\leftarrow} \mfk(\overline{\mb{Q}}; \ell ^i )\overset{\sim}{\rightarrow}
ku^G(*; \ell ^i)
\end{equation}
where $\mfk(\overline{\mb{Q}}; \ell ^i) $ denotes the spectrum $\mfk(Spec(\overline{\mb{Q}}), G) \wedge M_{\ell^i}$, and $ku^{G}(*; \ell ^i)$
denotes mod-$\ell ^i$ $G$-equivariant connective complex $K$-theory.  It is  understood that it is the colimit over the equivariant $K$-theories of the finite quotients defining $G$.  We now have the following diagram.  
 \begin{center} \vspace{-1cm}
\begin{equation} 
\begin{tikzcd}[row sep=2.5em, column sep=6em]
\Sigma ^2 \mfk(k; \ell ^i) \arrow[d,"\times \beta"]
& \Sigma ^2\mfk(\overline{\mb{Q}}; \ell ^i) \arrow[d,"\times \beta"] \arrow[r,"\sim"] \arrow[l, "\sim" above] &\Sigma ^2 ku^G(*; \ell ^i)\arrow[d,"\times \beta"]
 \\
 \mfk(k; \ell ^i) \arrow[d]
&   \mfk(\overline{\mb{Q}}; \ell ^i) \arrow[r,"\sim"] \arrow[l,"\sim" above] \arrow[d] & ku^G(*; \ell ^i) \arrow[d]\\
\mc{C}  {(k; \ell ^i)} & \mc{C} (\overline{\mb{Q}};\ell ^i)\arrow[l,"\sim" above] \arrow[r,"\sim"] & \mb{H}(\ell ^i)
\end{tikzcd}
\end{equation}
\end{center}
where the right hand $\times \beta $ is multiplication by the usual Bott map in connective complex $K$-theory.  The top two rows of maps are all equivalences by Suslin rigidity.  The spectra 
$\mc{C} {(k; \ell ^i)} $ and $\mc{C}(\overline{\mb{Q}};\ell ^i)$ are defined to be the cofibers of the maps $\times \beta:\Sigma ^2 \mfk(k; \ell ^i) \rightarrow \mfk(k; \ell ^i) $ and $\times \beta: \Sigma ^2 \mfk(\overline{\mb{Q}}; \ell ^i) \rightarrow \mfk(\overline{\mb{Q}}; \ell ^i)$, respectively.  It follows that the third row of maps are also equivalences.  

This discussion now shows that $\mc{C}(k; \ell ^i ) $ is an  Eilenberg-MacLane spectrum for the $R[G]$-module $R[G]/\ell ^i R[G]$.  It follows from Item \ref{derss} of Proposition \ref{properties} that the homotopy groups of $\mc{C} (k; \ell ^i)^{\wedge}_{\mb{F}_{\ell}}$ are isomorphic to the homotopy groups of the derived completion of the $R[G]$-module $R[G]/\ell ^i R[G]$ at the homomorphism $R[G] \rightarrow \mb{F}_{\ell }$. 

We now form the homotopy inverse limit of the diagrams (\ref{definediagram}). It follows from Lemma \ref{Moorestuff} that 
$$ R[G]^{\wedge}_{\mb{F}_{\ell}} \cong \underset{i }{holim}(R[G]/\ell ^i R[G])^{\wedge}_{\mb{F}_{\ell}}
$$
On the other hand, Lemma \ref{Moorestuff} also shows that
$$\mfk(k)^{\wedge}_{\mb{F}_{\ell}} \cong \underset{i}{holim} \hspace{.05cm}\mfk(k; \ell ^i)^{\wedge}_{\mb{F}_{\ell}} $$
and we conclude that the cofibre of $\beta: \Sigma ^2 \mfk (k)^{\wedge}_{\mb{F}_{\ell}} \rightarrow \mfk (k)^{\wedge}_{\mb{F}_{\ell}}$ is equivalent to the spectrum $R[G]^{\wedge}_{\mb{F}_{\ell}}$.  

We now turn to $\mfk(\overline{F})$.  Corollary \ref{finalcriterion} shows that we have an equivalence of $\bK(Spec(k))^{\wedge}_{\mb{F}_{\ell}}$-modules  
$$\mfk(Spec(F))^{\wedge}_{\bK(Spec(k)) \rightarrow \mb{F}_{\ell}} \rightarrow \mfk(\overline{F},G)^{\wedge} _{\bK\rightarrow \mb{F}_{\ell}} $$
This gives rise to a square
\vspace{-1cm}
\begin{center}
\begin{equation} \label{definediagramone}
\begin{tikzcd}[row sep=2.5em, column sep=4em]
\Sigma ^2 \mfk (Spec(F))^{\wedge}_{\bK(Spec(k))\rightarrow \mb{F}_{\ell} } \arrow[r,"\times \beta"] \arrow[d, "\sim"] & \mfk (Spec(F)) ^{\wedge} _{\bK(Spec(k))\rightarrow \mb{F}_{\ell}} \arrow[d,"\sim"]
 \\
\Sigma ^2 \mfk (E)^{\wedge}_{\bK \rightarrow \mb{F}_{\ell}}
\arrow[r, "\times \beta"] & \mfk (E)^{\wedge}_{\bK \rightarrow \mb{F}_{\ell}} 
\end{tikzcd}
\end{equation}
\end{center}
Our definition of the element $\beta$ is easily seen to give the usual Bott element by comparison with the continuous case $ku$.  The (proved) Beilinson-Lichtenbaum conjecture, as summarized in Theorem 5.2 of \cite{repassembly}, shows that the cofibre of the map in the top row of (\ref{definediagramone}) above has homotopy groups isomorphic to the $\ell$-adically completed Milnor $K$-groups of the field, which completes the proof.    
\end{Proof}

 \end{document}